\documentclass{article}

\usepackage[utf8]{inputenc}
\usepackage[T1]{fontenc}
\usepackage{geometry}
\usepackage{amsmath,amssymb,amsfonts,amsthm}
\usepackage{indentfirst}
\usepackage{graphicx}
\usepackage{xcolor}
\usepackage{verbatim}
\usepackage{mathrsfs}
\usepackage{mathbbol}
\usepackage{tikz}
\usetikzlibrary{patterns}
\usetikzlibrary{arrows}
\usepackage[english]{babel}
\usepackage{parskip}
\usetikzlibrary{shapes}
\usepackage{hyperref}
\usepackage{standalone}
\usepackage{authblk}
\usepackage{mathrsfs}

\usepackage{authblk} 
\usepackage[margin=16pt,font=small,labelfont=bf]{caption} 
\usepackage{float}
\usepackage{enumitem}

\newtheorem{thm}{Theorem}[section]
\newtheorem{proposition}[thm]{Proposition}
\newtheorem{corollary}[thm]{Corollary}
\newtheorem{lemma}[thm]{Lemma} 
\newtheorem{definition}[thm] {Definition}

\newtheorem{claim}[thm]{Claim}
\numberwithin{equation}{section} 
\newtheorem{remark}[thm]{Remark}

\newtheorem{notation}[thm]{Notation}

\definecolor{grigio}{RGB}{232,232,232}

 \newenvironment{theorem}{
 \begin{thm}}{\end{thm}}

 \newcommand{\be}[1]{\begin{equation}\label{#1}}
 \newcommand{\ee}{\end{equation}}

 \newcommand{\bl}[1]{\begin{lemma}\label{#1}}
 \newcommand{\el}{\end{lemma}}

 \newcommand{\br}[1]{\begin{remark}\label{#1}}
 \newcommand{\er}{\end{remark}}

 \newcommand{\bt}[1]{\begin{theorem}\label{#1}}
 \newcommand{\et}{\end{theorem}}

 \newcommand{\bd}[1]{\begin{definition}\label{#1}}
 \newcommand{\ed}{\end{definition}}

 \newcommand{\bcl}[1]{\begin{claim}\label{#1}}
 \newcommand{\ecl}{\end{claim}}

 \newcommand{\bp}[1]{\begin{proposition}\label{#1}}
 \newcommand{\ep}{\end{proposition}}

 \newcommand{\bc}[1]{\begin{corollary}\label{#1}}
 \newcommand{\ec}{\end{corollary}}

 \newcommand{\bn}[1]{\begin{notation}\label{#1}}
 \newcommand{\en}{\end{notation}}

  \newcommand{\bpr}{\begin{proof}}
  \newcommand{\epr}{\end{proof}}

 \newenvironment{proof*}{\noindent {\it Proof\,}}{\hspace{\fill}\qed\medskip}

 \newcommand{\bi}{\begin{itemize}}
 \newcommand{\ei}{\end{itemize}}

 \newcommand{\ben}{\begin{enumerate}}
 \newcommand{\een}{\end{enumerate}}

\def\sm{{\cX^m}}

 \def \ba {\begin{array}}
 \def \ea {\end{array}}

 \def \R {{\mathbb R}}
 
 \def \N {{\mathbb N}}
 \def \P {{\mathbb P}}
 \def \E {{\mathbb E}}

 \def \ra {\rightarrow}

 \def \cS {{\cal S}}
 \def \cE {{\cal E}}
 \def \cF {{\cal F}}
 \def \cA {{\cal A}}

\def \cI {{\cal I}}

\def \cK{{\cal K}}
 
 \def \cC {{\cal C}}
 \def \cD {{\cal D}}
 \def \cX {{\cal X}}

 \def \cW {{\cal W}}

 \def \cV {{\cal V}}

 \def \G {{\Gamma}}
 \def \L {{\Lambda}}

 \def \b {{\beta}}
 
 \def \D {{\Delta}}
 \def \r {{\rho}}

 \def \h {{\eta}}
 \def \s {{\sigma}}
 \def \z {{\zeta}}
 
 \def \t {{\tau}}
 \def \o {{\omega}}
 \def \d {{\delta}}

\def \x {{\xi}}

\def\GK{\G^{\text{K-Hex}}}

\newcommand{\bhexagon}{\mathord{\raisebox{0.6pt}{\tikz{\node[draw,scale=.65,regular polygon, regular polygon sides=6,fill=black](){};}}}}

\newcommand{\hexagon}{\mathord{\raisebox{0.6pt}{\tikz{\node[draw,scale=.65,regular polygon, regular polygon sides=6](){};}}}}

\def\pieno{{\bhexagon}}
\def\vuoto{{\hexagon}}

\begin{document}

\title{Metastability for Kawasaki dynamics \\ on the hexagonal lattice}

	\date{}
	
	\author[
	{}\hspace{0.5pt}\protect\hyperlink{hyp:email1}{1},\protect\hyperlink{hyp:affil1}{a,b}
	]
	{\protect\hypertarget{hyp:author1}{Simone Baldassarri}}
	
	\author[
	{}\hspace{0.5pt}\protect\hyperlink{hyp:email2}{2},\protect\hyperlink{hyp:affil1}{c}
	]
	{\protect\hypertarget{hyp:author2}{Vanessa Jacquier}}

	\affil[ ]{\centering
		\small\parbox{365pt}{\centering
			\parbox{5pt}{\textsuperscript{\protect\hypertarget{hyp:affil2}{a}}}Dipartimento di Matematica e Informatica ``Ulisse Dini", Universit\`{a} degli Studi di Firenze, Viale Morgagni 67/a, 50134 Firenze, Italy.
		}
	}
	
	\affil[ ]{\centering
		\small\parbox{365pt}{\centering
			\parbox{5pt}{\textsuperscript{\protect\hypertarget{hyp:affil2}{b}}}Aix-Marseille Universit\'e, CNRS, Centrale Marseille, I2M UMR CNRS 7373, 39, rue Joliot Curie, 13453 Marseille Cedex 13, France.
		}
	}
	
	\affil[ ]{\centering
		\small\parbox{365pt}{\centering
			\parbox{5pt}{\textsuperscript{\protect\hypertarget{hyp:affil2}{c}}}Scuola Normale Superiore, Piazza dei Cavalieri 7,
              56126 Pisa, Italy.
		}
	}
	
	\affil[ ]{\centering
		\small\parbox{365pt}{\centering
			\parbox{5pt}{\textsuperscript{\protect\hypertarget{hyp:email1}{1}}}\texttt{\footnotesize\href{mailto:simone.baldassarri@unifi.it}{simone.baldassarri@unifi.it}},
			\parbox{5pt}{\textsuperscript{\protect\hypertarget{hyp:email2}{2}}}\texttt{\footnotesize\href{mailto:vanessa.jacquier@sns.it}{vanessa.jacquier@sns.it}}
		}
	}

	\maketitle
	
	\vspace{-0.7cm}

 	\begin{center}
 		{\it Both authors are deeply grateful to their former group leader \\ Francesca Romana Nardi for her guidance and friendship, and for \\ encouraging the collaboration that led to this manuscript. \\ This work is dedicated to her memory.}		
 	\end{center}
	
\begin{abstract}
In this paper we analyze the metastable behavior for the Ising model that evolves under Kawasaki dynamics on the hexagonal lattice $\mathbb{H}^2$ in the limit of vanishing temperature. Let $\Lambda\subset\mathbb{H}^2$ a finite set which we assume to be arbitrarily large. Particles perform simple exclusion on $\Lambda$, but when they occupy neighboring sites they feel a binding energy $-U<0$. Along each bond touching the boundary of $\Lambda$ from the outside to the inside, particles are created with rate $\rho=e^{-\Delta\beta}$, while along each bond from the inside to the outside, particles are annihilated with rate 1, where $\beta$ is the inverse temperature and $\Delta>0$ is an activity parameter. For the choice $\Delta\in{(U,\frac{3}{2}U)}$ we prove that the empty (resp.\ full) hexagon is the unique metastable (resp.\ stable) state. 
We determine the asymptotic properties of the transition time from the metastable to the stable state and we give a description of the critical configurations. We show how not only their size but also their shape varies depending on the thermodynamical parameters. Moreover, we emphasize the role that the specific lattice plays in the analysis of the metastable Kawasaki dynamics by comparing the different behavior of this system with the corresponding system on the square lattice.

		\medskip
		\noindent
		{\it AMS} 2020 {\it subject classifications.} 60J10; 60K35; 82C20; 82C22; 82C26 
		
		
		\medskip
		\noindent
		{\it Acknowledgments.} This work was supported by the ``Gruppo Nazionale per l'Analisi Matematica e le loro Applicazioni" (GNAMPA-INdAM).

	\end{abstract}

\newpage

\section{Introduction}
Metastability is a dynamical phenomenon that occurs when a physical, chemical or biological system moves, under a stochastic dynamics, between different regions of its state space on different time scales. On short time scales the system is in {\it quasi-equilibrium} within a single region, while on long time scales it undergoes a sudden transition between {\it quasi-equilibria} within different regions. This transition is called {\it metastability} or {\it metastable behavior}. For very low temperature dynamics this phenomenon is characterized by the tendency of the system to remain trapped for extremely long time in a state (the metastable state) different from the stable states, until it performs a sudden transition to the stable state at some random time. The main issues that are typically investigated for a system that exhibits a metastable behavior are the following. The first is the analysis of the {\it transition time} from a metastable to a stable state. The second issue is the study of the so-called {\it critical configurations}, namely those configurations that are visited by the system with high probability during the transition. The last issue is the characterization of the {\it tube of typical trajectories} followed by the system with high probability during the crossover from the metastable to the stable states. 
This investigation has been carried over, in the literature, by using mainly two different approaches: the {\it pathwise} (see \cite{{CGOV},{OS1995},{OS1996},{OV}}) and the {\it potential theoretic} (see \cite{BEGK,BH}). Our results are obtained by leveraging on a modern version of the pathwise approach, which can be found in \cite{{CN2013},{CNS2015},{MNOS},{NZB},{Z},{Z2}}. The pathwise approach was applied in finite volume at low temperature in \cite{{AJNT2022},{BGK},{BGN},{BGNneg2021},{BGNpos2021},{CGOV},{CL},{NZ},{NS1992}} for single-spin-flip Glauber dynamics and in \cite{{CN2003},{CNS2008},{CNS2008bis}} for parallel dynamics, and we refer to \cite{CJS2022} for a review of the results for serial and parallel dynamics. The potential theoretic approach was applied to models at finite volume and at low temperature in \cite{{BJN},{BHN},{BM},{HNT2011},{HNT2012},{NS2012}}. The more involved infinite volume limit at low temperature or vanishing magnetic field was studied in \cite{{BHS},{CM2013},{GHNOS},{GMV},{GN},{HOS},{SS1998}} for Ising-like models under single-spin-flip Glauber and Kawasaki dynamics. More recent approaches are developed in \cite{{BL2010},{BL2015},{BG2016},{BGM2020}}.

In this paper we consider the metastable behavior of the two-dimensional isotropic Ising lattice gas at very low temperature and low density that evolves according to Kawasaki dynamics on the hexagonal lattice. Kawasaki dynamics is a discrete time Markov chain defined by the Metropolis algorithm with transition probabilities given in \eqref{defkaw}. Let $\b>0$ be the inverse temperature and let $\L\subset\mathbb{H}^2$ be a finite set such that its interior $\L^-$ is an hexagon (see Section \ref{model} for more details) with open boundary conditions. Particles live and evolve in a conservative way inside $\L$, but when they occupy neighboring sites they feel a binding energy $-U$. Along each bond touching the boundary of $\L$ from the outside to the inside, particles are created with rate $\rho=e^{-\D\b}$, while along each bond from the inside to the outside, particles are annihilated with rate 1, where $\D>0$ is an activity parameter. Thus, the boundary of $\L$ plays the role of an infinite gas reservoir with density $\rho$. We fix the parameters $U$ and $\D$ such that $\D\in(U,\frac{3}{2}U)$, that corresponds to the metastable regime. We will prove in Theorem \ref{thm:metastable_state} that the empty (resp.\ full) configuration is the unique metastable (resp.\ stable) state. We consider the asymptotic regime corresponding to finite volume $\L$ in the limit of large inverse temperature $\b$. We investigate how the system {\it nucleates}, i.e., how it reaches $\pieno$ (hexagon $\L^-$ full of particles) starting from $\vuoto$ (empty hexagon $\L^-$).

The main motivation of this paper is the following. From a physical point of view, the last two issues of metastability, namely the characterization of the critical configurations and the tube of typical trajectories, are the most relevant, because they provide a geometrical description of the evolution of the system. To this end, in this paper we investigate how the underlying lattice strongly affects the dynamical properties of the system. The choice of the hexagonal lattice comes from a recent study done for this model evolving under Glauber dynamics in \cite{{AJNT2022},{KS2022}}, because it has been shown how a certain class of parallel dynamics (shaken dynamics in \cite{{ADSST2019},{ADSST2019bis}}) on the square lattice induces a collection of parallel dynamics on a family of Ising models on the hexagonal lattice with non-isotropic interaction where the spins in each of the two partitions are alternatively updated. 

The goal of the paper is to investigate the critical configurations and the tunnelling time between $\vuoto$ and $\pieno$ for this model. To this end, in Section \ref{mainresults} we will give our main results: in Theorem \ref{thm:metastable_state} we identify the metastable and stable states. In Theorem \ref{thm:transition_time} we prove a convergence in probability, expectation and law for the transition time, answering the first issue introduced above. In Theorem \ref{prop:recurrence_property} we prove that the system reaches with high probability either the state $\vuoto$ or $\pieno$ in a time shorter than $e^{\beta(V^*+ \epsilon)}$, uniformly in the starting configuration for any $\epsilon> 0$, where $V^*=\Delta+U$. In other words, the dynamics speeded up by a factor of order $e^{\beta V^*}$ reaches with high probability $\{\vuoto, \pieno\}$. In Theorem \ref{thm:gate} we provide a characterization of a gate for the transition, namely a set of configurations which will be crossed with probability tending to one in the limit $\b\ra\infty$, answering the second issue of metastability. We emphasize that this result reflects how the underlying lattice is decisive for the dynamics of the system. One could be tempted to simply adapt the critical configurations for the same model on the square lattice to the hexagonal lattice, for example by replacing the rectangular shape with an hexagonal one, but this conjecture is false. Indeed, we will prove that for this model there exist two different sizes for the critical droplets depending on the value of the fractional part of the ratio $(\D-U)/(3U-2\D)$. This situation occurs also for the model evolving under Glauber dynamics considered in \cite{AJNT2022}, but we want to stress that its characterization is very different. Indeed, the main difference between Kawasaki and Glauber dynamics is that the former conserves the number of particles and therefore the structure of the gates is much richer. In particular, for Glauber dynamics there is a unique {\it minimal gate}, i.e., a gate minimal by inclusion, but for Kawasaki dynamics their characterization is not trivial and therefore much more interesting to derive. The geometrical description of the minimal gates is out of the scope of the present paper, but we encourage the reader to inspect the differences between Theorem \ref{thm:gate} and \cite[Theorem 2.13]{AJNT2022} for having in mind the different nature of the gate for the transition for these two different dynamics. We refer to Section \ref{comparison} for a detailed comparison between this model and the ones evolving under Kawasaki dynamics on the square lattice, in which we also emphasize that there are many shapes for the critical droplets according to the several kind of motions that can take place. Finally, in Theorem \ref{thm:subsup} we prove that, with probability tending to one, configurations with some hexagonal shape are subcritical, in the sense that they shrink to $\vuoto$ before reaching $\pieno$, or are supercritical, in the sense that they grow to $\pieno$ before reaching $\vuoto$. This result is a first step for the geometrical description of the tube of typical trajectories, namely the third issue of metastability.

\subsection{Comparison with Kawasaki dynamics on the square lattice}\label{comparison}
In this Section we make a comparison between the model we consider in this paper and other models evolving under Kawasaki dynamics on the square lattice in order to emphasize the different behavior of the system dependending on the geometry of the lattice. Indeed, this is the main motivation of the paper. There are many papers regarding the Ising lattice gas evolving under Kawasaki dynamics on the two and three-dimensional square lattice. For instance, in \cite{HOS} (resp.\ \cite{HNOS}) the isotropic version of this model is investigated in two (resp.\ three) dimensions by giving results concerning the asymptotics of the transition time and an intrinsic description of a gate. For two dimensions, in \cite{GOS} the complete description of the tube of typical paths is given and we also refer to \cite{{HNT2011},{HNT2012}} for the study of the model with two types of particles. Concerning the anisotropic version of this model, the weakly (resp.\ strongly) anisotropic case was first studied in \cite{NOS} (resp.\ \cite{BN2021strong}). In all these papers, an incomplete geometrical description of a gate for the transition from the empty box (metastable state) to the full box (stable state) is given. These results were sharpened in \cite{BHN} (resp.\ \cite{{BN2022strong},{BN2022weak}}) for the isotropic (resp.\ anisotropic) model both via the use of the potential theory both for a detailed geometrical characterization of the critical droplets. Indeed, a particular feature of Kawasaki dynamics is that in the metastable regime particles move along the border of the droplet more rapidly than they arrive from the boundary of the domain. The locally conservative dynamics and this movement of particles give a regularization effect, but we want to stress that the particular shape of the hexagonal lattice induces an increment of these regularizing motions in such a way new mechanisms of entering the critical configurations set appear, see Remarks \ref{treniniinterni} and \ref{remark:ingresso} for more details. This is a first crucial difference between the two isotropic models. Indeed, on the square lattice a new mechanism to enter the gate appears only in the strongly anisotropic setting, see \cite{{BN2021strong},{BN2022strong}}. For the weakly anisotropic and isotropic models there is a unique way to enter the gate: a rectangular shape with a single protuberance is reached and then a free particle enters from the boundary of the box, see \cite{{BN2022weak},{NOS}} for more details. On the square lattice, before the entrance of the free particle it is possible that particles move only along the border of the cluster, while on the hexagonal lattice this phenomenon can also appear for particles in an internal region of the cluster, see Figure \ref{trenino1}(a)-(b) for an example of the first and last configuration obtained in such a way. As a consequence, in this case the complete geometrical characterization is hard to obtain, and is left as a future research direction.
The reason we observe this very different behavior rests on the specific structure of the underlying lattice. Indeed, on the hexagonal lattice, when a particle that does not belong to the border of a cluster moves, if it attaches to a protuberance then the energy increases by $U$ (2 bonds are broken and one is created when the moving particle attaches to the protuberance), while this is false on the square lattice. Indeed, in that case the energy increases by $2U$ (3 bonds are broken and one is created when the moving particle attaches to the protuberance). This difference turns out to be crucial when the dynamics is close to critical configurations. This phenomenon can be also found in the different metastable regime for this model with respect to the one on the square lattice. This is peculiar of Kawasaki dynamics, indeed for Glauber dynamics this does not happen, see \cite[Condition 2.6]{AJNT2022}. We give an intuition of why this happens. For the two dimensional isotropic model the metastable regime corresponds to $\D\in(U,2U)$. Indeed, in this scenario single particles attached to one side of a droplet typically detach before the arrival of the next particle (because $e^{U\beta}\ll e^{\D\b}$), while bars of two or more particles typically do not detach (because $e^{\D\beta}\ll e^{2U\b}$). A similar interpretation can be derived for the analogous conditions which arise in the two-dimensional anisotropic cases and the three-dimensional isotropic case. For the hexagonal lattice the situation is different, indeed the metastable regime corresponds to $\D\in(U,\frac{3}{2}U)$. Clearly, the condition $\D>U$ has the same interpretation given above. But the condition $\D<2U$ is not enough as in the square lattice. Indeed, for more than one particle attached to an hexagonal shape is possible to detach a single particle alternatively at cost $U$ and $2U$ and therefore the required upper bound on $\D$ can be viewed as an average of these two costs. This particular behavior is also responsible for the particular shape of the critical droplets, which present two different protuberances and not only one as in the square lattice case. As it will be clear throughout the paper, we come to the conclusion that the geometry of the lattice significantly influences the behavior of the system subject to Kawasaki dynamics and this makes it very interesting to study.

\subsection{Outline of the paper}
The outline of the paper is as follows. In Section \ref{sec2} we define our model and the Kawasaki dynamics. Moreover, we give some model-independent and model-dependent definitions in order to state our main results. In Section \ref{stablevel} we give the proof of the theorems concerning the asymptotic behavior of the transition time and the characterization of the critical configurations after identifying the maximal stability level. This is done by providing an upper and lower bound via a reference path and by using the isoperimetric inequality respectively. Finally, in Section \ref{recurrencepropertyproof} we prove the recurrence property to the set $\{\vuoto,\pieno\}$ which allows use to identify the metastable and stable states for the system.

\section{Model and results}\label{sec2}

\subsection{Definition of the model}\label{model}
\begin{figure}
	\centering
	\includegraphics[scale=0.5]{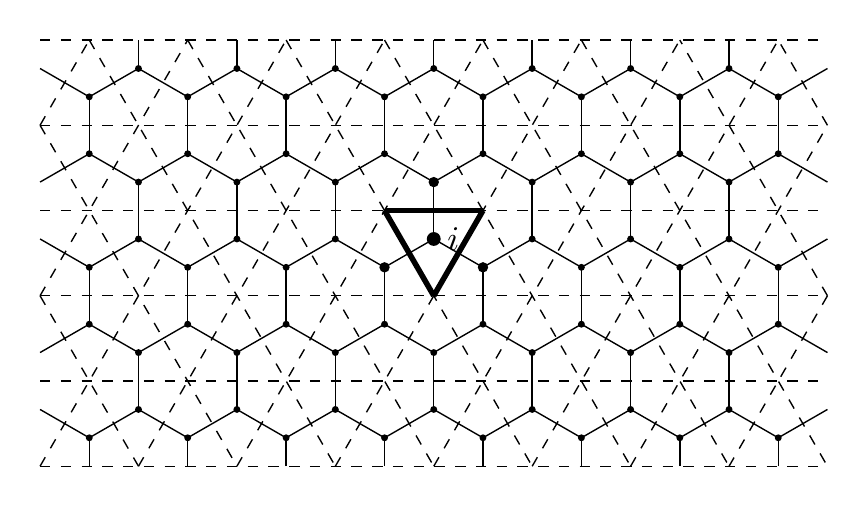}
	\hspace{0.05cm}
	\includestandalone[scale=0.3, mode=image|tex]{lambda}
	\caption{On the left-hand side we highlight in black the sites $j$ such that $d(i,j)\neq 0$ in our model. On the right-hand side we depict the set $\L$ with a straight line, while we depict the hexagon $\L^-$ with a dashed line.}
	\label{neighbors}
\end{figure}

Consider the discrete hexagonal lattice $\mathbb{H}^2$ embedded in $\mathbb{R}^2$ and let $\mathbb{T}^2$ be its dual, so that $\mathbb{T}^2$ is the triangular lattice. Two sites of the discrete hexagonal lattice are said to be \emph{nearest neighbors} when they share an edge of the lattice (see Figure \ref{neighbors} on the left-hand side). We consider an hexagon in $\mathbb{H}^2$ with radius $L$ and we define $\L\subset \mathbb{H}^2$ as the union between this hexagon and all the sites, that are not in the hexagon, with lattice distance one from the hexagon. Let
\begin{align}
	& \partial^- \Lambda := \{x \in \Lambda \, | \, \exists \, y \not \in \Lambda: |y - x| = 1\}
\end{align}
the internal boundary of $\Lambda$, and we put
\begin{align}
	& \Lambda^- := \Lambda \setminus \partial^-\Lambda.
\end{align}
With this choice of $\L$ we deduce that $\L^-$ is an hexagon with radius $L$ (see Figure \ref{neighbors} on the right-hand side). Note that $\L^-$ contains $6L^2$ sites. To each site $x\in\Lambda$ we associate an occupation variable $\h(x)\in \{0,1\}$ and we denote by $\h\in\mathcal{X}=\{{0,1}\}^{\Lambda}$ a lattice gas configuration. If the variable at site $x$ is zero, we say that the site is \emph{empty}, otherwise we say that it is \emph{occupied by a particle}.
On the \emph{configuration space} $\mathcal{X}$ we consider the \emph{Hamiltonian function} $H: \mathcal{X} \longrightarrow \mathbb{R}$ defined as
\begin{equation}\label{def:hamiltonian}
	H(\h):=-U\sum_{ \{x,y\}\in\Lambda^{*,-}} \h (x) \h (y) +\Delta \sum_{x \in \Lambda} \h (x),
\end{equation} 
where
\begin{equation}
	\Lambda^{*,-}= \{ \{x,y\} \in \Lambda^- \times \Lambda^-: \, |x -y| = 1\}
\end{equation}
is the set of non-oriented bonds in $\Lambda^-$. The interaction consists of a binding energy $-U<0$ for each nearest-neighbor pair of particles in $\Lambda^-$. In addition, there is an
activation energy $\Delta>0$ for each particle in $\Lambda$. The grand-canonical Gibbs measure associated with the Hamiltonian \eqref{def:hamiltonian} is
\begin{equation}\label{def:gibbs_measure}
	\mu_{\beta}(\h):=\frac{e^{-\beta H(\h)}}{Z_{\beta}}, \quad \hbox{with } \h\in\cX,
\end{equation}
where $Z_{\beta}:=\sum_{\eta \in \mathcal{X}} e^{-\beta H(\eta)}$ is the \emph{partition function} and $\beta:=\frac{1}{T} >0$ is the \emph{inverse temperature}.

\subsection{Local Kawasaki dynamics}\label{local_K_dynamics}

Next we define Kawasaki dynamics on $\L$ with boundary conditions that mimic the effect of an infinite gas reservoir outside $\L$ with density $ \r = e^{-\D\b}.$ Let $b=(x \to y)$ be an oriented bond, i.e., an {\it ordered} pair of nearest neighbour sites, and define
\begin{equation}\label{Loutindef}
	\begin{array}{lll}
		\partial^* \L^{out} &:=& \{b=(x \to y): x\in\partial^- \L,
		y\not\in\L\},\\
		\partial^* \L^{in}  &:=& \{b=(x \to y): x\not\in
		\L, y\in\partial^-\L\},\\
		\L^{*, orie} &:=& \{b=(x \to y): x,y\in\L\},
	\end{array}
\end{equation}
\noindent 
and put $ \bar\L^{*, orie}:=\partial^* \L ^{out}\cup
\partial^* \L ^{in}\cup\L^{*,\;orie}$.
Two configurations $  \sigma,
\sigma'\in {\cal X}$ with $ \sigma\ne \sigma'$ are said to be {\it
	communicating configurations} if there exists a bond
$b\in  \bar\L^{*,orie}$ such that $ \s' = T_b \s$, where $T_b   \sigma$ is the configuration obtained from $ \sigma$ in any of these ways:
\begin{itemize}
	\item
	For $b=(x \to y)\in\L^{*,\;orie}$, $T_b \sigma$ denotes the
	configuration obtained from $ \sigma$ by interchanging particles
	along $b$:
	\begin{equation}\label{Tint}
		T_b \sigma(z) =
		\left\{\begin{array}{ll}
			\sigma(z) &\mbox{if } z \ne x,y,\\
			\sigma(x) &\mbox{if } z = y,\\
			\sigma(y) &\mbox{if } z = x.
		\end{array}
		\right.
	\end{equation}
	\item
	For  $b=(x \to y)\in\partial^*\L^{out}$ we set:
	\begin{equation}\label{Texit}
		T_b \sigma(z) =
		\left\{\begin{array}{ll}
			\sigma(z) &\mbox{if } z \ne x,\\
			0     &\mbox{if } z = x.
		\end{array}
		\right.
	\end{equation}
	\noindent
	This describes the annihilation of a particle along the border.
	\item
	For  $b=(x  \to y)\in\partial^*\L^{in}$ we set:
	\begin{equation}\label{Tenter}
		T_b \sigma(z) =
		\left\{\begin{array}{ll}
			\sigma(z) &\mbox{if } z \ne y,\\
			1     &\mbox{if } z=y.
		\end{array}
		\right.
	\end{equation}
	\noindent
	This describes the creation of a particle along the border.
\end{itemize}

\noindent
The Kawasaki dynamics is  the discrete time Markov chain
$(\h_t)_{t\in \mathbb{N}}$ on state space $ {\cal X} $ given by
the following transition  probabilities: for  $  \h\not= \h'$:
\begin{equation}\label{defkaw}
P( \h,  \h'):=\left\{
\begin{array}{ll}
{ |\bar\L^{*,\;orie}|}^{-1} e^{-\b[H( \h') - H( \h)]_+}
&\text{if }  \exists \, b\in \bar\L ^{*, orie}: \h' =T_b \h,  \\
0   &\text{ otherwise, }  
\end{array} 
\right.
\end{equation}
where $[a]_+ =\max\{a,0\}$ and $P(\h,\h):=1-\sum_{\h'\neq\h}P(\h,\h')$. This describes a standard Metropolis dynamics with open boundary conditions: along each bond touching $\partial^-\L$ from the outside, particles are created with rate $\rho=e^{-\D\b}$ and are annihilated with rate 1, while inside $\L^-$ particles are conserved.
Note that an exchange of occupation
numbers $\h(x)$ for any $x$ inside $ \L\setminus  \L^-$
does not involve any change in energy.

\br{p1}
The stochastic dynamics defined by (\ref{defkaw}) is reversible with respect to the Gibbs measure in \eqref{def:gibbs_measure}.
\er

\subsection{Metastability: Static Heuristics}\label{section:euristica}
In this Section we present a heuristic discussion from a static point of view. We will consider the regime
\begin{equation}
	\Delta\in{\Big(U,\frac{3}{2}U\Big)}, \quad \beta\ra\infty,
\end{equation}
\noindent
which corresponds to the metastable behavior. Let us make a rough computation of the probability to see a regular hexagon of radius $r$ of occupied sites centered at the origin. We denote by $\mu^*$ the restricted ensemble, namely the Grand-canonical Gibbs measure defined in \eqref{def:gibbs_measure} restricted to a suitable subset of configurations, where all sufficiently large clusters are suppressed. Under this restricted ensemble we have
\begin{equation}
	\mu^*(\hbox{regular hexagon of radius } r)\approx \rho^{6r^2}e^{3U(3r^2-r)\beta},
\end{equation}
\noindent
since $\rho$ is close to the probability to find a particle at a given site and $-U$ is the binding energy between two particles at the neighboring sites, with $3(3r^2-r)$ the number of bonds for an hexagon with radius $r$. Writing $\rho=e^{-\Delta\beta}$ we obtain
\be{}
\mu^*(\hbox{regular hexagon of radius } r)\approx e^{-\beta[6r^2\Delta +3(r-3r^2)U]},
\ee
\noindent
where the exponent has a saddle point at
\be{}
\bar r=\frac{U}{2(3U-2\Delta)}.
\ee
This means that droplets with radius $r<\bar r$ have a tendency to shrink and droplets with radius $r\geq \bar r$ a tendency to grow. This would leave to the conclusion that $\bar r$ is the radius of the critical droplet. We will see in the sequel that the situation is more delicate (see \eqref{def:raggiocritico} for the precise definition of the critical radius $r^*$), indeed the dynamical mechanism for the transition between hexagonal droplets, which is not considered here, has an influence in establishing the tendency to grow or shrink. The choice $\Delta\in{(U,\frac{3}{2}U)}$ corresponds to $r^*\in{(1,\infty)}$, i.e., to a non-trivial critical droplet. The most interesting part of the metastable regime is $0<3U-2\Delta\ll U$, which corresponds to $r^*$ very large.

\subsection{Model-independent definitions and notations}\label{defmodind}
We will denote by $\P_{\h_0}$ the probability law of the Markov process  $(\h_t)_{t\geq 0}$ starting at $\h_0$ and by $\E_{\h_0}$ the corresponding expectation.

\noindent
{\bf 1. Paths and hitting times.}
\bi

\item
A {\it path\/} $\o$ is a sequence $\o=(\o_1,\dots,\o_k)$, with
$k\in\N$, $\o_i\in\cX$ and $P(\o_i,\o_{i+1})>0$ for $i=1,\dots,k-1$.
We write $\o\colon\;\h\to\h'$ to denote a path from $\h$ to $\h'$,
namely with $\o_1=\h,$ $\o_k=\h'$. Moreover, we denote by $\Theta(\h,\h')$ the set of all the paths connecting $\h$ and $\h'$. A set
$\cA\subset\cX$ with $|\cA|>1$ is {\it connected\/} if and only if for all
$\h,\h'\in\cA$ there exists a path $\o:\h\to\h'$ such that $\o_i\in\cA$
for all $i$. 

\item[$\bullet$]
Given a non-empty
set $\cA\subset\cX$, define the {\it first-hitting time of} $\cA$
as
\be{tempo}
\t_{\cA}:=\min \{t\geq 0:  \eta_t \in \cA \}.
\ee
\ei
\medskip
\noindent
\newpage
{\bf 2. Min-max and communication height}
\bi
\item Given a function $f:\cX\ra\R$ and a subset $\cA\subseteq\cX$, we denote by 
\be{defargmax}
\arg \hbox{max}_{\cA}f:=\{\sigma\in\cA: f(\sigma)=\max_{\z\in\cA}f(\z)\}
\ee
\noindent
the set of points where the maximum of $f$ in $\cA$ is reached. If $\o=(\o_1,...,\o_k)$ is a path, in the sequel we will write $\arg \max_{\o}H$ to indicate $\arg \max_{\cA}H$, with $\cA=\{\o_1,...,\o_k\}$ and $H$ the Hamiltonian.
\item 
The {\it bottom} $\cF(\cA)$ of a  non-empty
set $\cA\subset\cX$ is the
set of {\it global minima} of the Hamiltonian $H$ in  $\cA$:
\be{Fdef}
\cF(\cA):=\arg \hbox{min}_{\cA}H=\{\sigma\in\cA: H(\sigma)=\min_{\z\in\cA} H(\z)\}.
\ee
For a set $\cA\subset\cX$ such that all the configurations have the same energy, with an abuse of notation we denote this energy by $H(\cA)$.
\item 
The {\it communication height} between a pair $\sigma$, $\h\in\cX$ is
\be{}
\Phi(\sigma,\h):= \min_{\o:\sigma\ra\h}\max_{\z\in\o} H(\z).
\ee
\noindent
Given $\cA\subset\cX$, we define the {\it restricted communication height} between $\sigma,\h\in\cA$ as
\be{}
\Phi_{|\cA}(\sigma,\h):= \min_{\substack{\o:\sigma\ra\h\\ \o\subseteq\cA}}
\max_{\z\in\o} H(\z),
\ee
\noindent
where $(\o_1,...,\o_k)=\o\subseteq\cA$ means $\o_i\in\cA$ for every $i$.
\ei
\medskip
\noindent
{\bf 3. Stability level, stable and metastable states}
\bi
\item
We call
{\it stability level} of
a state $\sigma \in \cX$
the energy barrier
\be{stab}
V_{\sigma} :=
\Phi(\sigma,\cI_{\sigma}) - H(\sigma),
\ee
\noindent 
where $\cI_{\sigma}$ is the set of states with
energy below $H(\sigma)$:
\be{iz} 
\cI_{\sigma}:=\{\eta \in \cX \, | \, H(\eta)<H(\sigma)\}. 
\ee
\noindent 
We set $V_\sigma:=\infty$ if $\cI_\sigma$ is empty.
\item
We call {\it $V$-irreducible states}
the set of all states with stability level larger than  $V$:
\be{xv} 
\cX_V:=\{\h\in\cX \, | \, V_{\h}>V\}. 
\ee
\item
The set of {\it stable states} is the set of the global minima of
the Hamiltonian:
\be{st.st.}
\cX^s:=\cF(\cX). 
\ee
\item
The set of {\it metastable states} is given by
\be{st.metast.} 
\sm:=\{\sigma\in\cX \, | \,
V_{\sigma}=\max_{\eta\in\cX\setminus \cX^s}V_{\eta}\}. 
\ee
\noindent
We denote by $\G_m$ the stability level of the states in $\cX^m$.
\ei
\medskip
\noindent
{\bf 4. Optimal paths, saddles and gates}
\bi
\item 
We denote by $(\sigma\to\h)_{opt} $ the {\it set of optimal paths\/} as the set of all
paths from $\sigma$ to $\h$ realizing the min-max in $\cX$, i.e.,
\be{optpath}
(\sigma\to\h)_{opt}:=\{\o:\sigma\to\h\; \hbox{such that} \; \max_{\xi\in\o} H(\xi)=  \Phi(\sigma,\h) \}.
\ee
\item
The set of {\it minimal saddles\/} between
$\sigma,\h\in\cX$
is defined as
\be{minsad}
\cS(\sigma,\h):= \{\z\in\cX\colon\;\; \exists\o\in (\sigma\to\h)_{opt},
\ \o\ni\z \hbox{ such that } \max_{\xi\in\o} H(\xi)= H(\z)\}.
\ee
\item
Given a pair $\sigma,\h \in\cX$,
we say that $\cW\equiv\cW(\sigma,\h)$ is a {\it gate\/}
for the transition $\sigma\to\h$ if $\cW(\sigma,\h)\subseteq\cS(\sigma,\h)$
and $\o\cap\cW\neq\emptyset$ for all $\o\in (\sigma\to\h)_{opt}$. In words, a gate is a subset of $\cS(\sigma,\h)$ that is visited by all optimal paths.
\ei

\subsection{Model-dependent definitions}
We briefly give some model-dependent definitions and notations in order to state our main theorems. For the geometrical definitions see Section \ref{geomdef}. 
Recall that $\mathbb{T}^2$ is the dual of $\mathbb{H}^2$, i.e., $\mathbb{T}^2$ is the discrete triangular lattice embedded in
$\mathbb{R}^2$. 

\bi
\item[$\bullet$]
For $x\in\L^-$, let $ \hbox{nn}(x):=\{ y\in \L^-\colon\;|y-x| = 1\}$ be the set of nearest-neighbor sites of $x$ in $\L^-$ according to the lattice distance.
\item[$\bullet$]
A {\it free particle\/} in $\h\in\cX$ is a site $x$, with $\h(x)=1$, such that either $x\in\partial^-\L$, or $x\in\L^-$ and $\sum_{y\in nn(x)\cap\L^-}\h(y)$ $=0$. We denote by $\h_{free}$ the union of free particles in $\partial^-
\L$ and  free particles in $\L^-$. We denote by $n(\h)$ the number of free particles in $\h$.
\item[$\bullet$]
We denote by $\h_{cl}$ the clusterized
part of the occupied sites of $\h$:
\be{hcl}\h_{cl} :=\{x\in\L^-: \ \h(x)=1\}\setminus\h_{free}.
\ee
\item[$\bullet$]
We denote by $\h^{fp}$ the collection of configurations obtained by $\eta$ via the addition of a free particle anywhere in $\L$.
\item[$\bullet$]
We call \emph{triangular unit} or \emph{triangular face} an equilateral triangle of area one, whose center belongs to the discrete hexagonal lattice and whose vertices belong to its dual (see Figure~\ref{neighbors} on the left-hand side). Moreover, a set of two triangular units that share an edge is called \emph{elementary rhombus}.
\item[$\bullet$]
Given a configuration $\h \in \mathcal{X}$ we denote by $C(\h_{cl})$ its Peierls contour, that lives on the dual lattice and is the union of piecewise linear curves separating the empty triangular faces from the triangular faces with particles inside.
\item[$\bullet$]
Given a set $A\subset\mathbb{T}^2$, we define its area as the number of particles in $A$. We denote the area by $||A||$.
\item[$\bullet$]
The configuration space $\cX$ can be partitioned as
\be{}
\cX=\displaystyle\bigcup_{n}\cV_n,
\ee
\noindent
where $\cV_n:=\{\h\in\cX: \sum_{x\in\Lambda}\eta(x)=n\}$ is the set of configurations with $n$ particles, called the {\it $n$-manifold}.
\ei

\subsection{Main results}\label{mainresults}
In this section we present our main results for this model. Let
\be{def:vuoto}
\vuoto:=\{\h\in\cX: \h(x)=0 \ \forall \ x\in\L\}
\ee
\noindent
be the empty configuration. By \eqref{def:hamiltonian} and \eqref{def:vuoto} we have that $H(\vuoto)=0$. Let
\be{def:pieno}
\pieno:=\{\h\in\cX: \h(x)=1 \ \forall x\in\L^-, \ \h(x)=0 \ \forall x\in\L\setminus\L^-\}
\ee
\noindent
be the configuration that is full in $\L^-$ and empty in $\L\setminus\L^-$.
Define the critical radius $r^*$ as
 \begin{equation}\label{def:raggiocritico}
r^*:=\Big\lfloor\frac{U}{2(3U-2\Delta)}-\frac{1}{2}\Big\rfloor= \frac{U}{2(3U-2\Delta)}-\frac{1}{2}-\delta =\frac{\Delta-U}{3U-2\Delta}-\delta,
\end{equation}
\noindent
with $\d\in(0,1)$ the fractional part of $\frac{U}{2(3U-2\Delta)}-\frac{1}{2}$ fixed. We assume that $\frac{U}{2(3U-2\Delta)}-\frac{1}{2}$ not integer is made so to avoid strong degeneracy of the critical configurations. Similar assumptions are common in literature (see e.g., \cite{BM,CNS2008,CN2013}). We recall the assumption $3U-2\D\ll U$, in particular $3U-2\D\leq\frac{U}{100}$ is enough. In the following theorem, we will identify the stable and metastable states and we will show that for our model the energy barrier $\G_m$ is equal to
	\be{}\label{gammahex}
	\G^{\text{K-Hex}}:=
	\begin{cases}
 \Gamma^*_1 & \text{if } \d \in \big(0,\frac{1}{2}\big), \\
\Gamma^*_2 &\text{if } \d \in \big(\frac{1}{2},1\big).
	\end{cases}
	\ee
 where 
 \begin{equation*}
     \Gamma^*_1=-3(3(r^*)^2-r^*)U+6(r^*)^2\Delta+5(2r^*+1)\Delta-(15r^*+4)U+\D
     \end{equation*}
     and 
     \begin{equation*}
     \Gamma^*_2=-3(3(r^*+1)^2-(r^*+1))U+6(r^*+1)^2\Delta+(2r^*+3)\Delta-3(r^*+1)U+\D.
     \end{equation*}

\noindent
The value of $\G^{\text{K-Hex}}$ is obtained by computing the energy of the critical configurations. We will see that these configurations consist of a cluster having a shape that is close to a hexagon with radius $r^*$ and, in particular, we will compute the critical area to be
\be{def:criticalarea}
\ba{ll}
A_1^*=6(r^*)^2+10r^*+6 \qquad \qquad \qquad &\text{if } \d\in\big(0,\frac{1}{2}\big), \\
A_2^*=6(r^*+1)^2+2(r^*+1)+2  \qquad \qquad \qquad &\text{if } \d\in\big(\frac{1}{2},1\big).
\ea
\ee

\bt{thm:metastable_state}{\rm{(Identification of the metastable state).}}
Let $L>2r^*+3$, then $\cX^m=\{\vuoto\}$ and $\cX^s=\{\pieno\}$. Moreover, $\G_m=\Phi(\vuoto,\pieno)=\G^{\text{K-Hex}}$.
\et

The idea is to find an upper bound for $\Gamma_m$ by building a reference path and a lower bound using an isoperimetric inequality. Another of our goals is finding the asymptotic behavior as $\beta \to \infty$ of the transition time for the system started at the metastable state $\vuoto$. 
\bt{thm:transition_time} {\rm{(Asymptotic behavior of $\tau_{\pieno}$).}}
	For any $\epsilon>0$, we have
	\begin{equation}
		\lim_{\beta \to \infty} \mathbb{P}_{\vuoto}\Big(e^{\beta(\Gamma^{\text{K-Hex}}-\epsilon)}< \tau_{\pieno}<e^{\beta(\Gamma^{\text{K-Hex}}+\epsilon)}\Big) =1,
	\end{equation}
	\begin{equation}
		\lim_{\beta \to \infty}\dfrac{1}{\b} \log\mathbb{E}_{\vuoto} \tau_{\pieno}=\G^{\text{K-Hex}}.
	\end{equation}
	Moreover, letting $T_\b:=\inf\{n\geq1: \P_\vuoto(\t_\pieno\leq n)\geq1-e^{-1}\}$, we have
	\begin{equation}
		\lim_{\b\ra\infty}\P_{\vuoto}(\t_\pieno>tT_\b)=e^{-t}
	\end{equation}
	and
	\begin{equation}
		\lim_{\b\ra\infty}\dfrac{\E_{\vuoto}\t_\pieno}{T_\b}=1.
	\end{equation}
\et
We refer to Section \ref{gateproof} for the proof of Theorem \ref{thm:transition_time}. We say that a function $\beta \mapsto f(\beta)$ is super exponentially small (SES) if $$\lim_{\beta\to\infty}\frac{\log{f(\beta)}}{\beta}=-\infty.$$
With this notation we can state our first theorem concerning the recurrence of 
the system to either the state $\vuoto$ or $\pieno$.
\bt{prop:recurrence_property}{\rm{(Recurrence property).}} Let $V^*= \Delta+U$, we have \mbox{$\mathcal{X}_{V^*}\subseteq\{\vuoto, \pieno \}$} and for any $\epsilon>0$ and sufficiently large $\beta$, we have
	\begin{align}\label{eq:recurrence}
		\beta \mapsto \sup_{\sigma \in \mathcal{X}} \mathbb{P}_{\sigma}(\tau_{\mathcal{X}_{V^*}}> e^{\beta(V^*+\epsilon)}) \hbox{ is SES.}
	\end{align}
\et
Equation~\eqref{eq:recurrence} implies that the system reaches with high probability
either the state $\vuoto$ (which is a local minimizer of the Hamiltonian) or the ground state in a time shorter than 
$e^{\beta (V^*+\epsilon)}$, uniformly in the starting configuration $\sigma$ for any $\epsilon >0$. 
The proof of Theorem \ref{prop:recurrence_property} follows from Proposition \ref{prop:stability_lower} and \cite[Theorem 3.1]{MNOS} (see Section \ref{recurrencepropertyproof} for more details). 

\begin{figure}
	\centering
	\includegraphics[width=\textwidth]{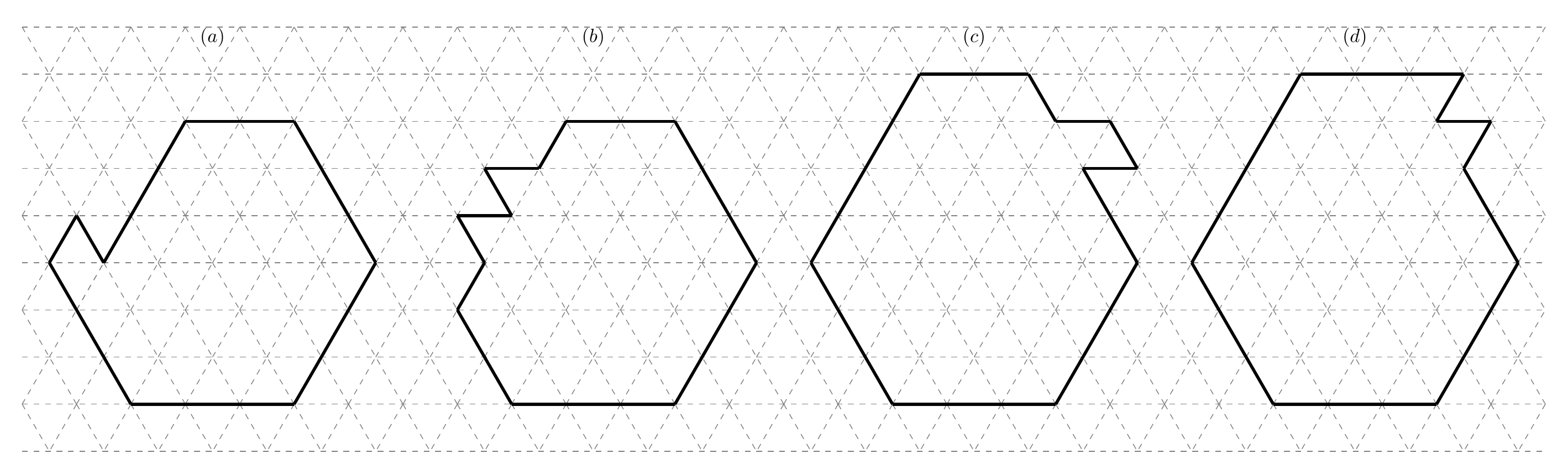}
	\caption{On the left there are two examples of configurations in $\mathcal{\tilde S}(A^*_1-1)$, $ \mathcal{\tilde D}(A^*_1-1)$ for $\delta \in (0,1/2)$, while on the right there are two examples of configurations in $\mathcal{\tilde S}(A^*_2-1)$, $\mathcal{\tilde D}(A^*_2-1)$ for $\delta \in (1/2,1)$.}
	\label{figselle}
\end{figure}

\begin{figure}[htb!]
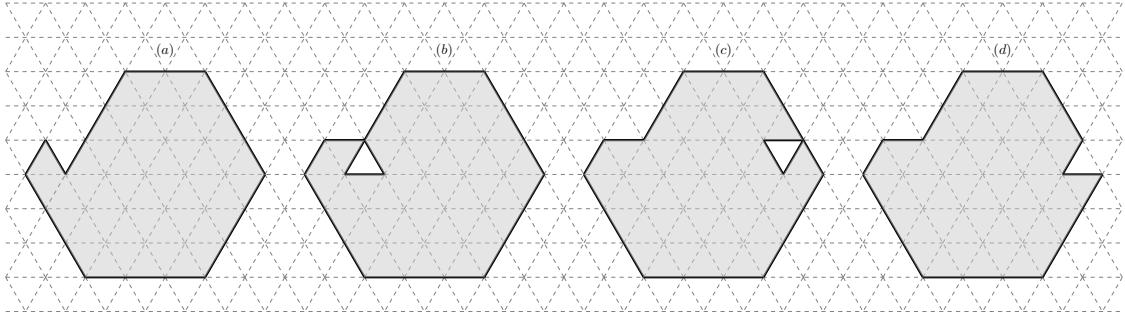

	\centering
	\includestandalone[width=\textwidth, mode=image|tex]{fig}
	\caption{In this figure we depict an example of motions of particles that belong to the internal part of a cluster at cost $U$. We represent the cluster in grey. Starting from the configuration represented in (a), by moving a particle towards the empty site, the energy increases by $U$ and the configuration that is obtained is the one represented in (b). From now on, the empty site moves at cost $0$ until the path reaches the configuration depicted in (c). Finally, the path reaches the configuration in (d) by lowering the energy by $U$, thus the starting and final configuration have the same energy.}
	\label{fignuova}
\end{figure}

In order to characterize the gate for the transition, we give an intuitive definition of the configurations denoted by $\tilde \cS(A_i^*-1)$ and $\tilde \cD(A_i^*-1)$ that play the role of protocritical configurations. In particular, configurations in $\tilde \cS(A_i^*-1)$ (resp.\ $\tilde \cD(A_i^*-1)$) have a unique cluster with area $A_i^*-1$ and shape as in Figure \ref{figselle}(a)(c)(resp.\ Figure \ref{figselle}(b)(d)). We refer the reader to Definitions \ref{standard} and \ref{defective} for a precise definition of these sets. Let
\begin{equation}\label{proto}
\begin{array}{ll}
	\cK(A_i^*-1):=\{&\h'\in\cV_{A_i^*-1}| \, \exists \, \omega=(\eta,\omega_1,...,\omega_n,\eta') \hbox{ such that } \\
 &\h\in \tilde \cS(A_i^*-1) \cup \tilde \cD(A_i^*-1), \, H(\h)=H(\h'), \\
 &n(\omega_j)=0 \, \forall \, j=1,...,n \hbox{ and } \Phi_{\cV_{A_i^*-1}}(\eta,\eta')\leq H(\h)+U\}
 \end{array}
\end{equation}
be the set of configurations obtained by a path starting from $\tilde \cS(A_i^*-1) \cup \tilde \cD(A_i^*-1)$ that conserves the number of particles and contains only configurations without free particles. Moreover, the energy along this path increases by $U$ at most and the starting and final configurations have the same energy. Note that the last condition in \eqref{proto} is the same as requiring that $\Phi_{\cV_{A_i^*-1}}(\eta,\eta')<\G^{K-Hex}$. The following theorem characterizes the gate for the transition from $\vuoto$ to $\pieno$. 
\bt{thm:gate}{\rm{(Gate for the transition).}}
	Given $\d\in(0,1)$ and $A_i^*\in\{A_1^*,A_2^*\}$ as in \eqref{def:criticalarea}, the set $\cC(A_i^*):= \cK(A_i^*-1)^{fp}$ is a gate for the transition from $\vuoto$ to $\pieno$.
\et

\br{treniniinterni}
Unlike what happens on the square lattice, on the hexagonal lattice much more ways to move particles at cost $U$ can take place. We want to stress this crucial property of the hexagonal lattice since it has a robust impact on the geometrical description of the gate. Indeed, for instance, concerning a configuration as in Figure \ref{fignuova}(a), note that it is possible to move a protuberance belonging to the elementary rhombus at cost $U$. The key fact is that these are not the unique possibilities, as occurs on the square lattice, indeed in this case it is possible to move also particles that belong to the internal part of a cluster. For example, it is possible to move towards the elementary rhombus an entire row of particles giving rise to a configuration with the same energy (see Figure \ref{fignuova} for the entire path). For this reason the geometrical characterization of the gate is much richer and more interesting than the one derived in \cite{BHN} for the square lattice. Moreover, these additional regularizing motions of particles lead to several mechanisms to enter the gate. We encourage the reader to inspect Remark \ref{remark:ingresso} for more details.
\er

We refer to Section \ref{gateproof} for the proof of Theorem \ref{thm:gate}.

In order to state the last result of this section, we need to introduce the set $\cE_{B_i}(r)$ that contains the configurations which have a unique cluster with a shape of \emph{quasi-regular hexagon}, that is a regular hexagon with $i$ bars attached clockwise. See figures \ref{bars} and \ref{regularhexagon} on the left-hand side and in the middle together with definitions \ref{hexagon}, \ref{def:bars}, \ref{qrhexagon} for more details.

\bt{thm:subsup}{\rm{(Subcritical and supercritical quasi-regular hexagons).}}
	Let $\cE_{B_i}^-(r)$ (resp.\ $\cE_{B_i}^+(r)$) be the set of configurations composed by a single quasi-regular hexagon contained in (resp.\ containing) $\cE_{B_i}(r)$.
	For $L>2r^*+3$, the following statements hold:
	\begin{itemize}
	\item[(i)] If $\d\in(0,\frac{1}{2})$, we have
	\begin{equation}
	\ba{ll}
	\hbox{if } \h\in\cE_{B_5}^{-} (r^*) &{ \implies } \displaystyle\lim_{\b\ra\infty}\P_{\h}(\t_\vuoto<\t_\pieno)=1, \\
	\hbox{if } \h\in\cE_{B_0}^{+} (r^*+1) &{ \implies } \displaystyle\lim_{\b\ra\infty}\P_{\h}(\t_\pieno<\t_\vuoto)=1.
	\ea
	\end{equation}
	\item[(ii)] If $\d\in(\frac{1}{2},1)$, we have
	\begin{equation}
	\ba{ll}
	\hbox{if } \h\in\cE_{B_1}^{-} (r^*+1) &{ \implies } \displaystyle\lim_{\b\ra\infty}\P_{\h}(\t_\vuoto<\t_\pieno)=1, \\
	\hbox{if } \h\in\cE_{B_2}^{+} (r^*+1) &{ \implies } \displaystyle\lim_{\b\ra\infty}\P_{\h}(\t_\pieno<\t_\vuoto)=1.
	\ea
	\end{equation}
	\end{itemize}
\et
In words, we characterize subcritical and supercritical quasi-regular hexagons, i.e., subcritical quasi-regular hexagons shrink to $\vuoto$, while supercritical quasi-regular hexagons grow to $\pieno$. We refer to Section \ref{subsup} for the proof of Theorem \ref{thm:subsup}.

\section{Identification of maximal stability level}\label{stablevel}

\subsection{Geometrical definitions}\label{geomdef}
Now we recall some geometrical definitions and properties about clusters and polyiamonds present in \cite{AJNT2022}.

\bd{def:polyiamond}
A \emph{polyiamond} $P \subset \mathbb{R}^2$ is a finite maximally connected union of three or more triangular units that share at least a side. 
\ed

Note that if two triangular units share only a point these are considered, by definition, two different polyiamonds. 

We define a new bijection that associates to each cluster a polyiamond with the same shape. This implies that to each cell without a particle, we associate an empty triangular unit.  

\bd{def:perimeter_of_polyiamond}
	The \emph{boundary} of a polyiamond $P$ is the collection of unit edges of the lattice $\mathbb{T}^2$ such that each edge separates a triangular unit belonging to $P$ from an empty triangular unit. The \emph{edge-perimeter} $p(P)$ of a polyamond $P$ is the cardinality of its boundary. 
\ed
In other words the perimeter is given by the number of interfaces on the discrete triangular lattice $(\mathbb{T}^2)$ between the sites inside the polyiamond and those outside. 
If not specified differently, we will refer to the edge-perimeter simply as perimeter. 
\bd{def:internal_perimeter}
	The \emph{external boundary} of a polyiamond consists of the connected components of the boundary such that for each edge there exists a hexagonal-path in $\mathbb{H}^2$ which connects this edge with the boundary of $\L$ without intersecting
	the polyiamond. The \emph{internal boundary} of a polyiamond consists of the connected components of the boundary that are not external. 
\ed

\begin{figure}
		\centering \includegraphics[scale=0.5]{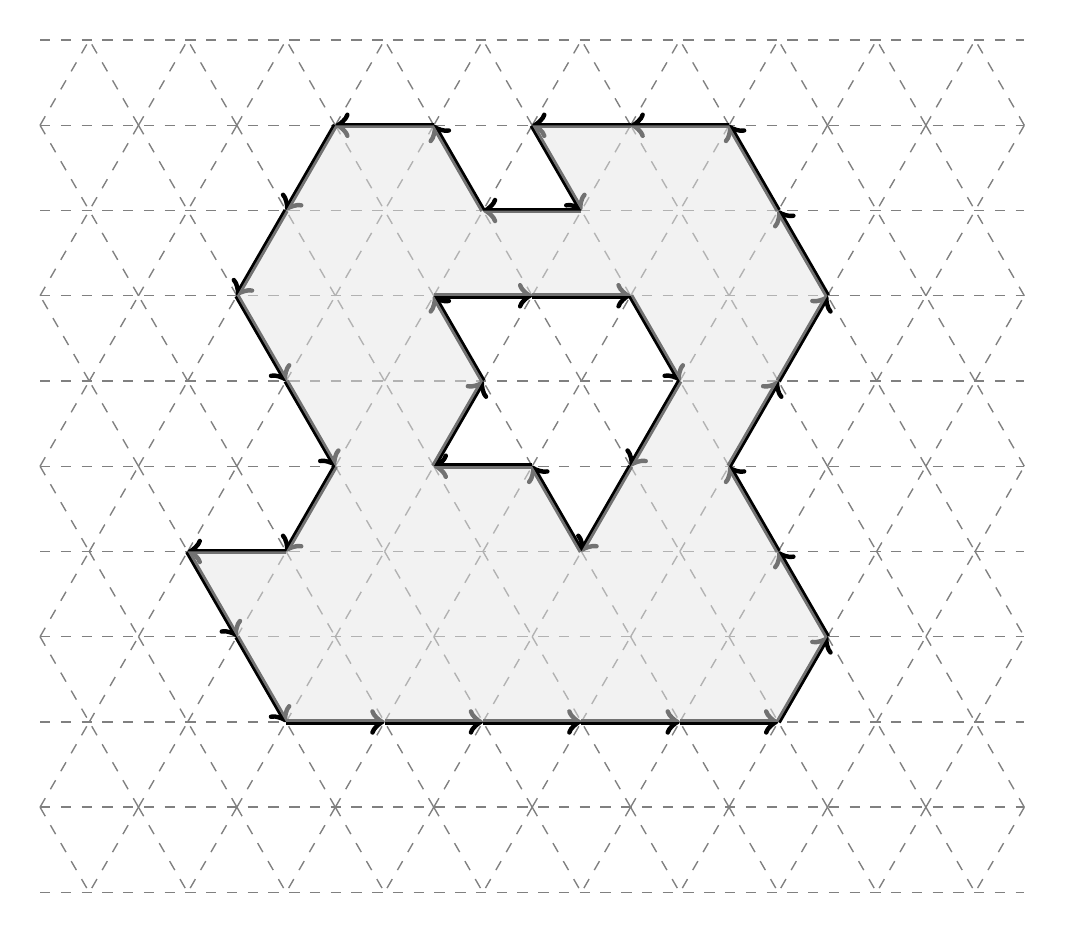} 
		\caption{An example of polyiamond with the external boundary oriented counter-clockwise and the internal boundary oriented clockwise.} \label{anglesfigure}
	\end{figure}

\bd{internalangle}
	Let us orient counter-clockwise the external boundary and clockwise the internal boundary. For each pair of oriented edges, the angle defined rotating counter-clockwise the second edge on the first edge is called \emph{internal angle} (see Figure \ref{anglesfigure}).
\ed
\bd{def:hole_in_polyiamond}
	A \emph{hole} of a polyiamond $P$ is a finite maximally connected component of empty triangular units surrounded by the internal boundary of $P$.  
\ed
We refer to holes consisting of a single empty triangle as \emph{elementary holes}.
\bd{}
	A polyiamond is \emph{regular} if it has only internal angles of $\pi$ and $\frac{2}{3}\pi$ and it has no holes. 
\ed
We note that a regular polyiamond has the shape of a hexagon.

	\begin{figure}
		\centering \includegraphics[scale=0.45]{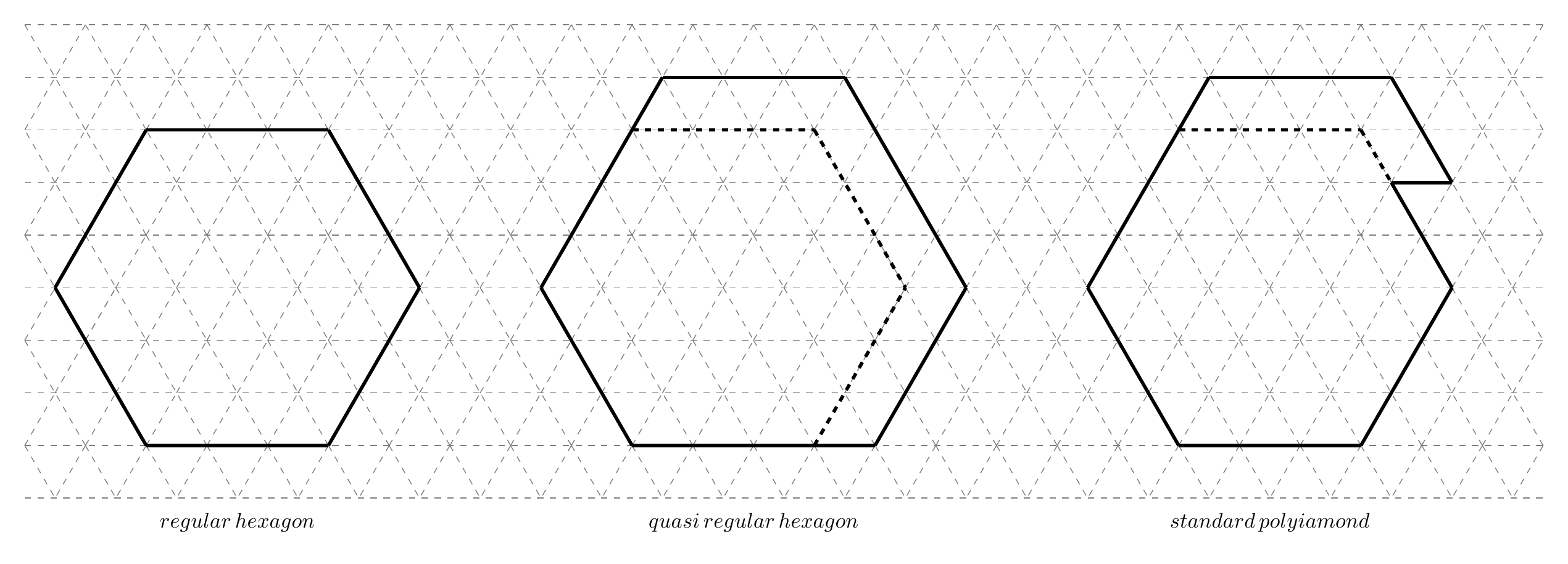} 
		\caption{Starting from the left: a regular hexagon with radius length 3, a quasi-regular hexagon with three bars attached clockwise from the top, a standard polyiamond with a complete bar attached on the top and an incomplete bar with cardinality 4 attached close to the top.}
		\label{regularhexagon}
	\end{figure}

\bd{hexagon}
A polyiamond is a \emph{regular hexagon} if it is a regular polyiamond with all equal sides. We denote by $E(r)$ the regular hexagon, where $r$ is its radius (see Figure \ref{regularhexagon} on the left-hand side).
\ed

\begin{figure}
	\centering \includegraphics[scale=0.5]{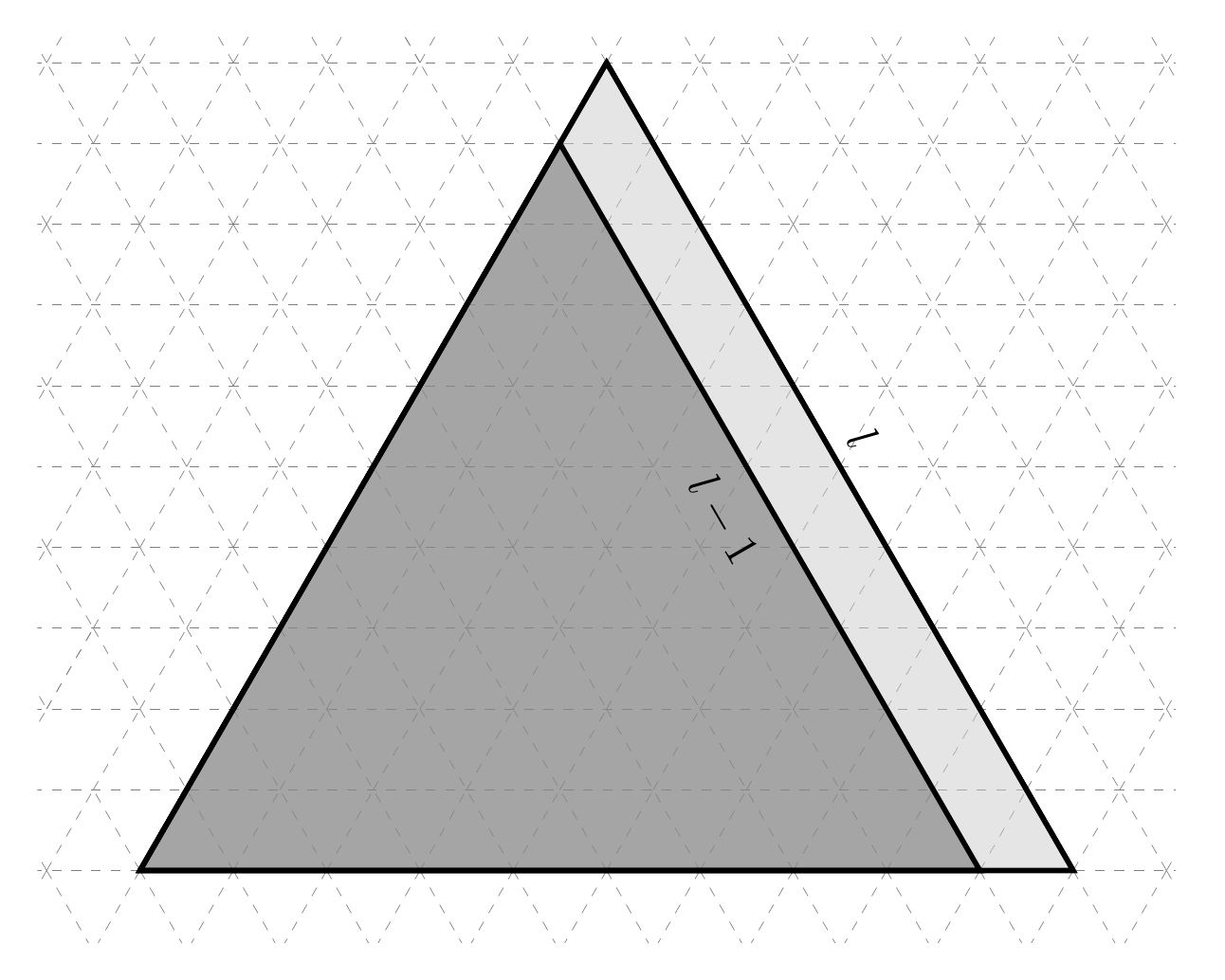} 
	\caption{The lightly shaded triangular units form a bar with larger base $l$ and
		it is obtained as difference between an equilateral triangle with
		side length $l$ and an equilateral triangle with side length $l-1$.}
	\label{fig:barra}
\end{figure}

\bd{def:bars}
A \emph{bar} $B(l)$ with larger base $l$ is a set of $||B(l)||=2l-1$ triangular units obtained as a difference between an equilateral triangle with side length $l$ and another equilateral triangle with side length $l-1$ (see Figure \ref{fig:barra}).
\ed

\begin{figure}
	\centering \includegraphics[width=\textwidth]{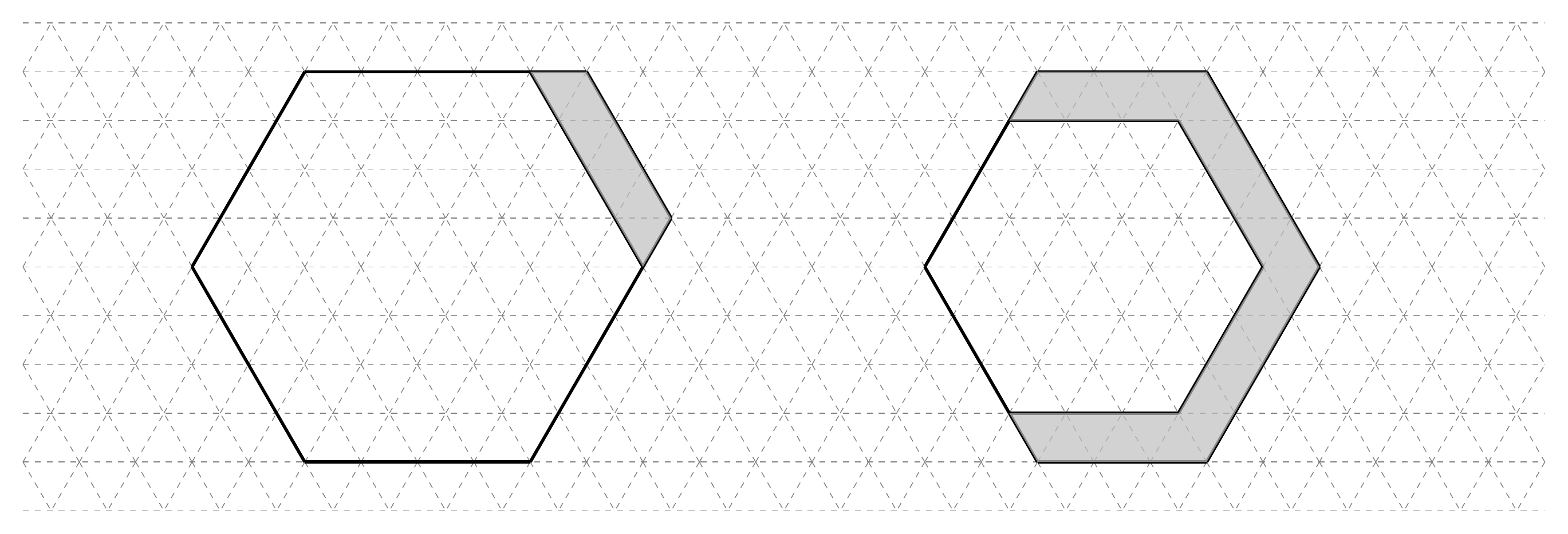} 
	\caption{On the left a quasi-regular hexagon $E_{B_1}(4)$. We observe that the cardinality of $B_1$ of $E_{B_1}(3)$ is $||B_1||=2r-1=5$. On the right a quasi-regular hexagon $E_{B_4}(3)$. We observe that the cardinality of $B_1$ of $E_{B_1}(3)$ is $||B_1||=2r-1=5$, while the cardinality of $B_i$ of $E_{B_4}(3)$ is $||B_i||=2r+1=7$ with $i=2,\ldots,4$.} \label{bars}
\end{figure}

\bd{qrhexagon}
	We denote by $E_{B_1}(r)$ the polyiamond obtained attaching a bar $B_1$ along its larger base $r$ to the regular hexagon (see Figure \ref{bars}). Analogously, we denote by $E_{B_i}(r)$ for $i=2,...,5$ the polyiamonds obtained attaching a bar $B_i$ along its larger base $r+1$ to $E_{B_{i-1}}(r)$.	Finally, we denote by $E_{B_6}(r)$ the polyiamond obtained attaching a bar $B_6$ along its larger base $r+2$ to $E_{B_5}(r)$. We call $E_{B_i}(r)$ a \emph{quasi-regular hexagon}, where $r$ is the radius of the regular hexagon $E(r)$ and $i\in \{1,...,6\}$ is the number of bars attached to it. 
\ed

Note that $E_{B_i}(r)$ is always contained in $E(r+1)$ and it is defined up to a rotation of $z \frac{\pi}{3}$ for $z \in \mathbb{Z}$.
Moreover $E(r) \equiv E_{B_0}(r)$ and $E(r+1) \equiv E_{B_6}(r)$. 

\bn{}
We denote by $\cE(r)$ the set of configurations $\h\in\cX$ such that $\h$ has a unique cluster with shape $E(r)$. We denote by $\cE_{B_i}(r)$ the set of configurations $\h\in\cX$ such that $\h$ has a unique cluster with shape $E_{B_i}(r)$.
\en

\begin{figure}[htb!] 
	\centering \includegraphics[scale=0.4]{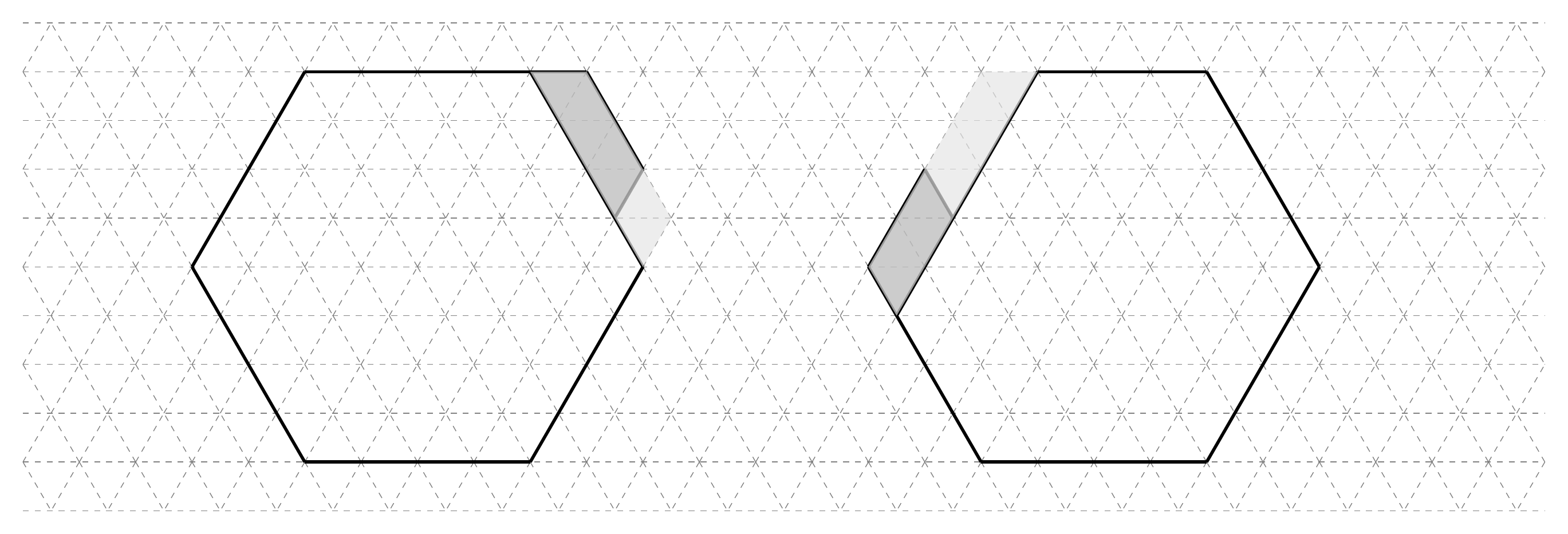} 
	\caption{On the left an incomplete bar with trapeze shape and cardinality $k=5$ attached to the regular hexagon $E(4)$. We observe that the cardinality of the bar containing the incomplete bar is $||B_1||=7>k$. On the right an incomplete bar with parallelogram shape and cardinality $k=4$, attached to the quasi-regular hexagon $E_{B_5}(3)$. We observe that the cardinality  of the bar containing the incomplete bar is $||B_6||=9>k$.} \label{incompletebar}
\end{figure}

\bd{}
An \emph{incomplete bar} of cardinality $k < 2l-1$ is a subset of a bar with larger base $l$ (see Figure \ref{incompletebar}). 
\ed

\bd{standard}
A \emph{standard polyiamond} of area $A$, denoted by $S(A)$, is a quasi-regular hexagon $E_{B_i}(r)$ with possibly an additional incomplete bar of cardinality $k$ attached clockwise, such that it is contained in $E_{B_{i+1}}(r)$. 
If $k=2$, we denote it by $\tilde S(A)$.
\ed

When we refer to a \emph{standard cluster} with area $A$, our meaning is that the cluster has the shape and the properties of a standard polyiamond $S(A)$.

Note that a standard polyiamond $S(A)$ is determined solely by its area $A$. We characterize $S(A)$ in terms of its radius $r$, the number $i$ of bars attached to the regular hexagon $E(r)$ to obtain $E_{B_i}(r)$ and the cardinality $k$ of the possible incomplete bar. Starting from the area $A$, these values can be computed by using \cite[algorithm 3.18]{AJNT2022}. 

\bd{defective}
A polyiamond consisting of a quasi-regular hexagon with two triangular units attached to one of its longest sides at triangular lattice distance 2 one from the other is called \emph{defective} and it is denoted by $\tilde D(A)$, where $A$ is the area.
\ed

\bn{}
We denote by $\tilde\cS(A)$ (resp.\ $\tilde\cD(A)$) the set of configurations $\h\in\cX$ such that $\h$ has a unique cluster with shape $\tilde S(A)$ (resp.\ $\tilde D(A)$). See (a)(c) (resp.\ (b)(d)) in Figure \ref{figselle} for examples of configurations in $\tilde\cS(A)$ (resp.\ $\tilde\cD(A)$). 
\en

\bd{corner}
We call a \emph{corner} of a standard polyiamond $P$ the pair of triangular faces of $P$ contained in the internal angle of $\frac{2}{3}\pi$.
\ed

\subsection{Reference path}\label{refpath} 
In this Section we construct our {\it reference path} $\o^*$, which is a sequence of configurations connecting $\vuoto$ and $\pieno$ such that the maximum of the energy along this path is $\G^{\text{K-Hex}}$. In particular, this path is composed by increasing clusters as close as possible to quasi-regular hexagonal shape. The idea of the construction of $\o^*$ is the following: we first construct a skeleton path $\{\bar\o_r\}_{r=0}^L$ given by a sequence of configurations that contain regular hexagon with radius $r$. Obviously $\bar\o_r$ is not a path in the sense that the transition from $\bar\o_r$ to $\bar\o_{r+1}$ can not be given in a single step of the dynamics, since $\bar\o_r$ and $\bar\o_{r+1}$ contain regular hexagons. Thus in order to obtain a path we have to interpolate each transition of the skeleton path. This is done in two steps. First, we introduce between $\bar\o_r$ and $\bar\o_{r+1}$ a sequence of configurations $\tilde\o_{r}^1,...,\tilde\o_{r}^{i_r}$ given by $\bar\o_r$ plus a bar, i.e., a quasi-regular hexagon. Again, these configurations are given by a single increasing droplet. Finally, we introduce a second interpolation to obtain a path $\o^*$ from the sequence of configurations $\tilde\o_{r}^i$. Its construction goes as follows. Between every couple of consecutive configurations in $\tilde\o$, for which the cluster is increased by one particle, a sequence of configurations with a new particle is inserted. In particular, the new particle is initially created at the boundary of $\L$ and then brought to the correct site via consecutive moves of this free particle.

{\bf Skeleton $\bar\o$:} Let us construct a sequence of configurations that contain regular hexagons $\bar\o=\{\bar\o_r\}$, with $r=0,...,L$, such that $\bar\o_0=\vuoto,..., \bar\o_L=\pieno$ and $\bar\o_r\subset\bar\o_{r+1}$. Starting from the origin, given $r=1,...,L$ let $\bar\o_r$ the regular hexagon with radius $r$, i.e., $\bar\o_r\in\cE(r)$. 

{\bf First interpolation $\tilde\o$:} From $\bar\o_0$ to $\bar\o_1$, we define the path $\tilde\o_{0}^i$ such that $\tilde\o_{0}^0=\bar\o_0$ and insert between $\bar\o_0$ and $\bar\o_1$ a sequence of configurations $\{\tilde\o_{0}^{i}\}_{i=0}^{6}$, which correspond to the creation of a hexagon of radius one obtained by adding sequentially particles clockwise. Given a choice for $\bar\o_r$, with $r<L$, we can construct the path $\tilde\o_{r}^{i}$ such that $\tilde\o_{r}^0=\bar\o_r$ and insert between $\bar\o_r$ and $\bar\o_{r+1}$ a sequence of configurations $\{\tilde\o_{r}^{i}\}_{i=0}^{i_r}$ as follows. Starting from a configuration in $\cE(r)$, add consecutive triangular units to the regular hexagon until we obtain a configuration in $\cE_{B_1}(r)$. Next we fill the bar on the top right adding consecutive triangular units until we obtain a configuration in $\cE_{B_2}(r)$. We go on in the same way adding bars clockwise, until we obtain configurations in $\cE_{B_3}(r),...,\cE_{B_6}(r)\equiv\cE(r+1)$. 

{\bf Second interpolation $\o^*$:} For any pair of configurations $(\tilde\o_{r}^{i},\tilde\o_{r}^{i+1})$ such that $||\tilde\o_{r}^{i}||<||\tilde\o_{r}^{i+1}||$, by construction of the path $\tilde\o_{r}^{i}$ the particles are created along the external boundary of the clusters, except for the first particle that is at the origin. So there exist $x_1,...,x_{j_i}$ a connected chain of nearest-neighbor empty sites of $\tilde\omega_{r}^{i}$ such that $x_1\in\partial^-\L$ and $x_{j_i}$ is the site where is located the additional particle in $\tilde\omega_{r}^{i+1}$. Define
\be{}
\hat\o_{r}^{i,0}=\tilde\o_{r}^{i}, \qquad \hat\o_{r}^{i,j_i}=\tilde\o_{r}^{i+1}, \quad r=0,...,L.
\ee
\noindent
Insert between each pair $(\tilde\o_{r}^{i},\tilde\o_{r}^{i+1})$ a sequence of configurations $\hat\o_{r}^{i,j}$, with $j=1,...,j_i-1$, where the free particle is moving from $x_1\in\partial^-\L$ to the cluster until it reaches the position $x_{j_i}$. Otherwise, for any pair of configurations $(\tilde\o_{r}^{i},\tilde\o_{r}^{i+1})$ such that $||\tilde\o_{r}^{i}||=||\tilde\o_{r}^{i+1}||$, we define $\hat\o_{r}^{i,0}=\tilde\o_{r}^{i}$ and $\hat\o_{r}^{i+1,0}=\tilde\o_{s}^{i+1}$. This concludes the definition of the reference path.

With an abuse of notation we denote a configuration in $\cE_{B_i}(r)$ by $\cE_{B_i}(r)$.
\bp{prop:quasiregular}
The maximum of the energy in $\o^*$ between two consecutive quasi-regular hexagons $\Phi_{\o^*}(\cE_{B_i}(r),\cE_{B_{i+1}}(r))$ for every $i=0,...,5$ is achieved in the standard polyiamond obtained adding to $\cE_{B_i}(r)$ an elementary rhombus along the longest side and a free particle.
\ep

\bpr
	Let $A^{(n)}$ be the area obtained after adding $n$ triangular units to the area of the quasi-regular hexagon in $\cE_{B_i}(r)$, where $n=0,...,||B_{i+1}||$. Note that $A^{(n)}$ is the area of the standard polyiamond $\cS(A^{(n)})$. We observe that $S(A^{(0)})=E_{B_i}(r)$ and $S(A^{(||B_{i+1}||)})=E_{B_{i+1}}(r)$. Since the maximum of the energy is obtained after adding a free particle, we obtain
	\begin{equation}
		H(\cS(A^{(n)})^{fp})-H(\cS(A^{(n-1)})^{fp})=
		\begin{cases}
			\Delta-U &\text{if } n=1, \\
			\Delta-U &\text{if } n \text{ is even}, \\
			\Delta-2U &\text{if } n\neq1 \text{ is odd}.
		\end{cases}
	\end{equation}
	Therefore we deduce that
	\begin{equation} 
		H(\cS(A^{(n)})^{fp})-H(\cE_{B_i}(r)^{fp})=
		\begin{cases}
			U-n\Big(\frac{3}{2}U-\Delta\Big) &\text{if } n \text{ is even}, \\
			\dfrac{U}{2}-n\Big(\frac{3}{2}U-\Delta\Big) &\text{if } n \text{ is odd}.
		\end{cases}
	\end{equation}
	Since the r.h.s. of the last equation decreases with $n$, due to the fact that $\D<\frac{3}{2}U$, in both the odd and even case, it is immediate to check that the maximum is attained for $n=2$, namely in $\cS(A^{(2)})^{fp}$.
\epr

\bp{togliere}
The maximum of the energy in $\o^*$ between two consecutive quasi-regular hexagons $\Phi_{\o^*}(\cE_{B_i}(r),\cE_{B_{i-1}}(r))$ for every $i=1,...,6$ is achieved in the standard polyiamond obtained removing counter-clockwise from $\cE_{B_i}(r)$ a number of particles equals to $||B_i||-3$ and detaching the $(||B_i||-2)$-th particle from ${B_i}$.
\ep
\begin{proof}
	Let $A^{(n)}$ be the area obtained after adding $n$ triangular units to the area of the quasi-regular hexagon $E_{B_{i-1}}(r)$, where $n=0,...,||B_i||$. Note that $\cS(A^{(n)})$ can be obtained either by removing $||B_i||-n$ triangular units from $\cE_{B_i}(r)$ or by adding and attaching $n$ triangular units to the quasi-regular hexagon in $\cE_{B_{i-1}}(r)$. We recall that {\it removing} a triangular unit means detaching it from the cluster and moving the free particle outside $\L$. Since the maximum of the energy is obtained after adding a free particle, we obtain
	\be{}
	H(\cS(A^{(n-1)})^{fp})-H(\cS(A^{(n)})^{fp})=
	\begin{cases}
		2U-\Delta &\text{if } n\neq1 \text{ is odd}, \\
		U-\Delta &\text{if } n \text{ is even}, \\
		U-\Delta &\text{if } n=1.
	\end{cases}
	\ee
	Therefore we deduce that
	\be{}
	H(\cS(A^{(n)})^{fp})-H(\cE_{B_i}(r)^{fp})=
	\begin{cases}
		n\Big(\dfrac{3}{2}U-\Delta\Big)-U &\text{if } n \text{ is even}, \\
		n\Big(\dfrac{3}{2}U-\Delta\Big)-\dfrac{U}{2} &\text{if } n \text{ is odd}.
	\end{cases}
	\ee
	Since the r.\ h.\ s.\ of the last equation increases with $n$, due to the fact that $\Delta<\frac{3}{2}U$, in both the odd and even case, it is immediate to check that the maximum is attained by removing $||B_i||-3$ triangular units from $\cE_{B_i}(r)$ and detaching another triangular unit from $B_i$. Therefore we obtain a configuration in $\cS(A^{(2)})^{fp}$.
\end{proof}

Recalling \eqref{def:raggiocritico}, from now on the strategy is to divide the reference path $\o^*$ into three regions depending on $r$:
\bi
\item the region $r\leq r^*$ will be considered in Proposition \ref{prop:sottocritico};
\item the region $r=r^*+1$ will be considered in Proposition \ref{prop:critico};
\item the region $r\geq r^*+2$ will be considered in Proposition \ref{prop:supercritico}
\ei

\bp{prop:sottocritico}
If $r\leq r^*$, then the communication height between two consecutive regular hexagons $\Phi_{\o^*}(\cE(r),\cE(r+1))$ along the path $\o^*$ is achieved in a configuration with a free particle and a standard cluster such that the number of its triangular units is $\tilde A=6r^2+10r+5$, that is $\Phi_{\o^*}(\cE(r),\cE(r+1))=\Phi_{\o^*}(\cE_{B_5}(r),\cE(r+1))=H(\cS(\tilde A))+\Delta$. Moreover, $\Phi_{\o^*}(\vuoto,\cE(r^*+1))=\Phi_{\o^*}(\cE(r^*),\cE(r^*+1))=H(\cS(A^*_1-1))+\Delta$ is achieved in a configuration with a free particle and a standard cluster $S(A^*_1-1)$, where $A^*_1=6(r^*)^2+10r^*+6$.
\ep

\begin{proof}
	Let $S(A)$ be a standard polyiamond with an incomplete bar of cardinality two. We obtain:
 {\small{	
 \begin{equation}\label{energie1}
	H(\cS(A))=
	\begin{cases}
		-3(3r^2-r)U+6r^2\Delta+2(\Delta-U) &\text{if } A=6r^2+2, \\
		-3(3r^2-r)U+6r^2\Delta+(2r+1)\Delta-3rU &\text{if } A=6r^2+2r+1, \\
		-3(3r^2-r)U+6r^2\Delta+2(2r+1)\Delta-(6r+1)U &\text{if } A=6r^2+4r+2, \\
		-3(3r^2-r)U+6r^2\Delta+3(2r+1)\Delta-(9r+2)U &\text{if } A=6r^2+6r+3, \\
		-3(3r^2-r)U+6r^2\Delta+4(2r+1)\Delta-(12r+3)U &\text{if } A=6r^2+8r+4, \\
		-3(3r^2-r)U+6r^2\Delta+5(2r+1)\Delta-(15r+4)U &\text{if } A=6r^2+10r+5. 
	\end{cases}
	\end{equation}
 }}
	\noindent
	We compare $\Phi_{\o^*}(\cE(r),\cE_{B_1}(r))=\Phi_{\o^*}(\cS(6r^2),\cS(6r^2+2r-1))$ with $\Phi_{\o^*}(\cE_{B_1}(r),\cE_{B_2}(r))= \Phi_{\o^*}(\cS(6r^2+2r-1),\cS(6r^2+4r))$. By Proposition \ref{prop:quasiregular} we have:
	\begin{equation}{}
	\ba{ll}
	\Phi_{\o^*}(\cE(r),\cE_{B_1}(r))&=H(\cS(6r^2+2))+\Delta, \\
	\Phi_{\o^*}(\cE_{B_1}(r),\cE_{B_2}(r))&=H(\cS(6r^2+2r+1))+\Delta.
	\ea
	\end{equation}
	\noindent
	By using \eqref{energie1}, we obtain that $\Phi_{\o^*}(\cE(r),\cE_{B_1}(r))\leq\Phi_{\o^*}(\cE_{B_1}(r),\cE_{B_2}(r))$ if and only if $r\leq\frac{2U-\D}{3U-2\D}=\frac{U}{2(3U-2\D)}+\frac{1}{2}$, which is true since we are assuming $r\leq r^*$ and $r^*\leq\frac{2U-\D}{3U-2\D}$ due to the condition $2\D< 3U$.
	
	We compare $\Phi_{\o^*}(\cE_{B_1}(r),\cE_{B_2(r)})=\Phi_{\o^*}(\cS(6r^2+2r-1),\cS(6r^2+4r))$ with $\Phi_{\o^*}(\cE_{B_2}(r),\cE_{B_3(r)})=\Phi_{\o^*}(\cS(6r^2+4r),\cS(6r^2+6r+1))$. By Proposition \ref{prop:quasiregular} we have:
	\be{}
	\ba{ll}
	\Phi_{\o^*}(\cE_{B_1}(r),\cE_{B_2}(r))&=H(\cS(6r^2+2r+1))+\D, \\
	\Phi_{\o^*}(\cE_{B_2}(r),\cE_{B_3}(r))&=H(\cS(6r^2+4r+2))+\D.
	\ea
	\ee
	By using \eqref{energie1}, we obtain that $\Phi_{\o^*}(\cE_{B_1}(r),\cE_{B_2}(r))\leq\Phi_{\o^*}(\cE_{B_2}(r),\cE_{B_3}(r))$ if and only if $r\leq\frac{\D-U}{3U-2\D}$, which is true since we are assuming $r\leq r^*$. 
	
	We compare $\Phi_{\o^*}(\cE_{B_2}(r),\cE_{B_3(r)})=\Phi_{\o^*}(\cS(6r^2+4r),\cS(6r^2+6r+1))$ with $\Phi_{\o^*}(\cE_{B_3}(r),\cE_{B_4(r)})=\Phi_{\o^*}(\cS(6r^2+6r+1),\cS(6r^2+8r+2))$. By Proposition \ref{prop:quasiregular} we have:
	\be{}
	\ba{ll}
	\Phi_{\o^*}(\cE_{B_2}(r),\cE_{B_3}(r))&=H(\cS(6r^2+4r+2))+\D, \\
	\Phi_{\o^*}(\cE_{B_3}(r),\cE_{B_4}(r))&=H(\cS(6r^2+6r+3))+\D.
	\ea
	\ee

	By using \eqref{energie1}, we obtain that $\Phi_{\o^*}(\cE_{B_2}(r),\cE_{B_3}(r))\leq\Phi_{\o^*}(\cE_{B_3}(r),\cE_{B_4}(r))$ if and only if $r\leq\frac{\D-U}{3U-2\D}$, which is true since we are assuming $r\leq r^*$.
	
	By performing similar computations, we obtain the following inequalities:
	\be{}
	\ba{ll}
	\Phi_{\o^*}(\cE_{B_3}(r),\cE_{B_4}(r))&\leq\Phi_{\o^*}(\cE_{B_4}(r),\cE_{B_5}(r)), \\
	\Phi_{\o^*}(\cE_{B_4}(r),\cE_{B_5}(r))&\leq\Phi_{\o^*}(\cE_{B_5}(r),\cE(r+1)).
	\ea
	\ee
	Thus the communication height between two consecutive regular hexagons along the path $\o^*$ is achieved in $\cS(6r^2+10r+5)^{fp}$, that is $\Phi_{\o^*}(\cE(r),\cE(r+1))=\Phi_{\o^*}(\cE_{B_5}(r),\cE(r+1))=H(\cS(\tilde A))+\D$, where $\tilde A=6r^2+10r+5$. The maximum of the function $H(\cS(\tilde A))+\D=-3(3r^2-r)U+6r^2\Delta+5(2r+1)\Delta-(15r+4)U+\D$ is obtained in $r=\frac{U}{2(3U-2\D)}-\frac{5}{6}$. However $r\in\N$ and $r\leq r^*$, therefore the maximum is attained in $r^*$ and $\Phi_{\o^*}(\vuoto,\cE(r+1))=\Phi_{\o^*}(\cE(r^*),\cE(r^*+1))=H(\cS(A^*_1-1))+\Delta$, where $A^*_1=6(r^*)^2+10r^*+6$.
\end{proof}

\bp{prop:critico}
If $r=r^*+1$, then the communication height $\Phi_{\o^*}(\cE(r^*+1),\cE(r^*+2))$ along the path $\o^*$ is achieved in a configuration with a free particle and a standard cluster $S(A^*_2-1)$, where $A^*_2=6(r^*+1)^2+2(r^*+1)+2$.
\ep

\begin{proof}
	Note that in this case $r=\Big\lfloor \frac{U}{2(3U-2\D)}+\frac{1}{2}\Big\rfloor$. We analyze $\Phi_{\o^*}(\cE(r^*+1),\cE(r^*+2))$ by using Proposition \ref{prop:quasiregular}. We compare the same communication height of Proposition \ref{prop:sottocritico}, obtaining the following inequalities since $r=r^*+1$:
	\be{}
	\Phi_{\o^*}(\cS(6r^2),\cS(6r^2+2r-1))<\Phi_{\o^*}(\cS(6r^2+2r-1),\cS(6r^2+4r)),
	\ee
	and
	\be{}
	\ba{ll}
	\Phi_{\o^*}(\cS(6r^2+2r-1),\cS(6r^2+4r))&>\Phi_{\o^*}(\cS(6r^2+4r),\cS(6r^2+6r+1))\\
	&>\Phi_{\o^*}(\cS(6r^2+6r+1),\cS(6r^2+8r+2))\\
	&>\Phi_{\o^*}(\cS(6r^2+8r+2),\cS(6r^2+10r+3))\\
	&>\Phi_{\o^*}(\cS(6r^2+10r+3),\cS(6r^2+12r+6)).
	\ea
	\ee
	Then the communication height along the path $\o^*$ between two consecutive regular hexagons with radius $r^*+1$ is $\Phi_{\o^*}(\cE(r^*+1),\cE(r^*+2))=\Phi_{\o^*}(\cS(6(r^*+1)^2+2(r^*+1)-1),\cS(6(r^*+1)^2+4(r^*+1)))$ and, by Proposition \ref{prop:quasiregular}, it is attained in $\cS(A^*_2-1)^{fp}$, with $A^*_2=6(r^*+1)^2+2(r^*+1)+2$. 
\end{proof}

\bp{prop:supercritico}
If $r\geq r^*+2$, then the communication height between two consecutive regular hexagons $\Phi_{\o^*}(\cE(r),\cE(r+1))$ along the path $\o^*$ is achieved in a configuration with a free particle and a standard cluster such that the number of its triangular units is $\tilde A=6r^2+2$, that is $\Phi_{\o^*}(\cE(r),\cE(r+1))=\Phi_{\o^*}(\cE(r),\cE_{B_1}(r))=H(\cS(\tilde A))+\D$. Moreover, $\Phi_{\o^*}(\cE(r^*+2),\pieno)=\Phi(\cE(r^*+2),\cE(r^*+3))=H(\cS(A^*_3-1))+\D$ is achieved in a configuration with a free particle and a standard cluster $S(A^*_3-1)$, where $A^*_3=6(r^*+2)^2+3$.
\ep

\begin{proof}
	We analyze $\Phi_{\o^*}(\cE(r),\cE(r+1))$ by using Proposition \ref{prop:quasiregular}. We compare the same communication height of Proposition \ref{prop:sottocritico}, obtaining the following inequalities since $r\geq r^*+2$:
	\be{}
	\ba{ll}
	\Phi_{\o^*}(\cS(6r^2),\cS(6r^2+2r-1))&\geq\Phi_{\o^*}(\cS(6r^2+2r-1),\cS(6r^2+4r)) \\
	&\geq\Phi_{\o^*}(\cS(6r^2+4r),\cS(6r^2+6r+1))\\
	&\geq\Phi_{\o^*}(\cS(6r^2+6r+1),\cS(6r^2+8r+2))\\
	&\geq\Phi_{\o^*}(\cS(6r^2+8r+2),\cS(6r^2+10r+3))\\
	&\geq\Phi_{\o^*}(\cS(6r^2+10r+3),\cS(6r^2+12r+6)).
	\ea
	\ee
	Thus the communication height between two consecutive regular hexagons along the path $\o^*$ is attained in $\cS(6r^2+2)^{fp}$, that is $\Phi_{\o^*}(\cE(r),\cE(r+1))=\Phi_{\o^*}(\cE(r),\cE_{B_1}(r))=H(\cS(\tilde A))+\D$, where $\tilde A=6r^2+2$. The maximum of the function $H(\cS(\tilde A))+\D=-3(3r^2-r)U+6r^2\D+2(\D-U)+\D$ is attained in $r=\frac{U}{2(3U-2\D)}$, but $r\in\N$ and $r\geq r^*+2$, so $\Phi_{\o^*}(\cE(r^*+2),\pieno)=\Phi_{\o^*}(\cE(r^*+2),\cE(r^*+3))=H(\cS(A^*_3-1))+\D$, where $A^*_3=6(r^*+2)^2+3$.
\end{proof}

\begin{figure}
	\centering \includegraphics[scale=0.5]{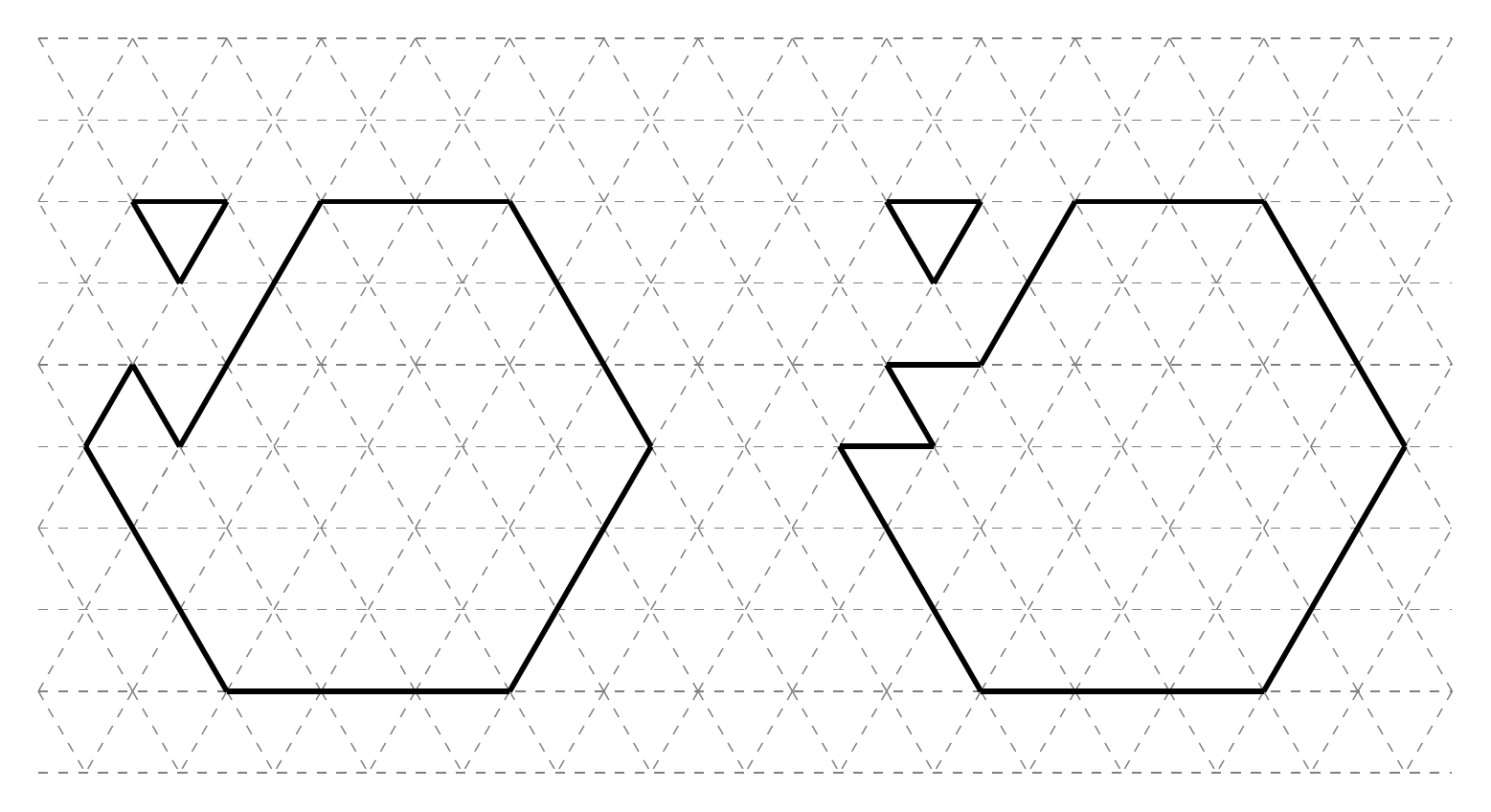} 
	\caption{Two standard clusters with critical area $A_1^*-1$ and a free particle for $0<\d<\frac{1}{2}$.} \label{selleKawasaki}
\end{figure}

\bp{prop:upperbound}
Let $\d\in(0,1)$ be such that $\frac{U}{2(3U-2\D)}-\frac{1}{2} -\d$ is an integer number. The maximum $\Phi_{\o^*}(\vuoto,\pieno)$ along the path $\o^*$ is attained in a configuration with a free particle and a standard cluster with area $A_i^*-1$ for $i\in\{1,2\}$ (see Figure \ref{selleKawasaki}) , where
\bi
\item[1)] $A_1^*=6(r^*)^2+10r^*+6$ if $0<\d<\frac{1}{2}$;
\item[2)] $A_2^*=6(r^*+1)^2+2(r^*+1)+2$ if $\frac{1}{2}<\d<1$.
\ei
\ep
\begin{proof}
	We compare $\Phi_{\o^*}(\vuoto,\cE(r^*+1))$, $\Phi_{\o^*}(\cE(r^*+1),\cE(r^*+2))$ and $\Phi_{\o^*}(\cE(r^*+2),\pieno)$. By Proposition \ref{prop:sottocritico} we have
	\be{eq1}
	\ba{ll}
	\Phi_{\o^*}(\vuoto,\cE(r^*+1))&=H(\cS(6(r^*)^2+10r^*+5))+\D \\
	&=-3(3(r^*)^2-r^*)U+6(r^*)^2\D+5(2r^*+1)\D-(15r^*+4)U+\D.
	\ea
	\ee
	By Proposition \ref{prop:critico} we have
	\be{eq2}
	\ba{ll}
	\Phi_{\o^*}(\cE(r^*+1),\cE(r^*+2))&=H(\cS(6(r^*+1)^2+2(r^*+1)+1))+\D \\
	&=-3(3(r^*+1)^2-(r^*+1))U+6(r^*+1)^2\D \\
	& \,\,\,\,\,\, +(2(r^*+1)+1)\D-3(r^*+1)U+\D.
	\ea
	\ee
	By Proposition \ref{prop:supercritico} we have
	\be{eq3}
	\ba{ll}
	\Phi_{\o^*}(\cE(r^*+2),\pieno)&=H(\cS(6(r^*+2)^2+2))+\D \\
	&=-3(3(r^*+2)^2-(r^*+2))U+6(r^*+2)^2\D+2(\D-U)+\D.
	\ea
	\ee
	By comparing equations \eqref{eq1},\eqref{eq2} and \eqref{eq3}, we obtain
	\be{}
	\ba{ll}
	\Phi_{\o^*}(\vuoto,\cE(r^*+1))&>\Phi_{\o^*}(\cE(r^*+2),\pieno),\\
	\Phi_{\o^*}(\cE(r^*+1),\cE(r^*+2))&>\Phi_{\o^*}(\cE(r^*+2),\pieno).
	\ea
	\ee
	Thus we deduce that $\Phi_{\o^*}(\cE(r^*+2),\pieno)$ cannot be the maximum. Moreover, we obtain
	\begin{align}
		\ba{ll}
		\Phi_{\o^*}(\vuoto,\cE(r^*+1))>\Phi_{\o^*}(\cE(r^*+1),\cE(r^*+2)) \quad &\text{if } 0<\d<\frac{1}{2}, \\
		\Phi_{\o^*}(\vuoto,\cE(r^*+1))<\Phi_{\o^*}(\cE(r^*+1),\cE(r^*+2)) \quad &\text{if } \frac{1}{2}<\d<1.
		\ea
	\end{align}
	and therefore the maximum $\Phi_{\o^*}(\vuoto,\pieno)=\Phi_{\o^*}(\vuoto,\cE(r^*+1))$ is achieved in a configuration $\cS(6(r^*)^2+10r^*+5)^{fp}$ if $\d \in (0,\frac{1}{2})$. Otherwise, if $\d\in(\frac{1}{2},1)$, then the maximum $\Phi_{\o^*}(\vuoto,\pieno)=\Phi_{\o^*}(\cE(r^*+1),\cE(r^*+2))$ is achieved in a configuration $\cS(6(r^*+1)^2+2(r^*+1)+1)^{fp}$.
\end{proof}

\bc{cor:upperG}
	Let $\GK$ as in \eqref{gammahex}. We have
	\be{}
	\Phi(\vuoto,\pieno)\leq\GK.
	\ee
\ec

\begin{proof}
By definition of communication height and the fact that $H(\vuoto)=0$, Proposition \ref{prop:upperbound} implies that
\be{}
\Phi(\vuoto,\pieno)\leq\max_{i}H(\o_i^*)=\GK
\ee
in the two cases $0<\d<\frac{1}{2}$ and $\frac{1}{2}<\d<1$.
\end{proof}

\subsection{Lower bound of maximal stability level}
In this Section we will find a lower bound for $\G^{\text{K-Hex}}$. In particular, we prove that $\Phi(\vuoto,\pieno)\geq\G^{\text{K-Hex}}$ separately for the case $\d \in \big(0,\frac{1}{2}\big)$ and $\d \in \big(\frac{1}{2},1\big)$. The proof comes in three steps, which are the contents of the three following lemmas. The last result of this section combines the upper and lower bound on $\Phi(\vuoto,\pieno)$ which we have found.

\bl{lemma:lower1}
The following statements hold:
\bi
\item[1.] If $\d \in \big(0,\frac{1}{2}\big)$, any $\o\in(\vuoto\ra\pieno)_{opt}$ must pass through the set $\mathcal{E}_{B_5}(r^*)$.
\item[2.] If $\d \in \big(\frac{1}{2},1\big)$, any $\o\in(\vuoto\ra\pieno)_{opt}$ must pass through the set $\mathcal{E}_{B_1}(r^*+1)$.
\ei
\el

\begin{proof}
	We analyze separately the two cases.
	\bi
	\item[1.] Let $\d \in \big(0,\frac{1}{2}\big)$ and $\tilde A=6(r^*)^2+10r^*+3$. Any path $\o:\vuoto\ra\pieno$ must cross the set $\cV_{\tilde A}$. By using \cite[Theorem 3.22]{AJNT2022} and \cite[Lemma 3.24]{AJNT2022} with $m=5$, in $\cV_{\tilde A}$ the unique (modulo translations and rotations) configuration of minimal perimeter and hence minimal energy is the standard polyiamond $S(\tilde A)$, which contains only the quasi-regular hexagon. Thus, the configuration $\cS(\tilde A)$ has energy
	\be{}
	\ba{ll}
	H(\cS(\tilde A))&=-3(3(r^*)^2-r^*)U+6(r^*)^2\Delta+5(2r^*+1)\Delta-(15r^*+1)U\\
	&=\G^{\text{K-Hex}}-3\D+2U.
	\ea
	\ee
	All the other configurations in $\cV_{\tilde A}$ have energy at least $\G^{\text{K-Hex}}-3\D+3U$. To increase the particle number starting from any such a configuration, we must create a particle at cost $\D$. But the resulting configuration would have energy $\G^{\text{K-Hex}}-2\D+3U$, which exceeds $\G^{\text{K-Hex}}$ due to the condition $2\D<3U$. Thus this would lead to a path exceeding the energy value $\G^{\text{K-Hex}}$ and therefore the path would not be optimal.
	\item[2.] Let $\d \in \big(\frac{1}{2},1\big)$ and $\tilde A=6(r^*+1)^2+2(r^*+1)-1$. By observing that \cite[Lemma 3.24]{AJNT2022} holds with $m=1$, we can argue as before.
	\ei
\end{proof}

\bl{lemma:lower2}
	The following statements hold:
	\bi
	\item[1.] If $\d \in \big(0,\frac{1}{2}\big)$, any $\o\in(\vuoto\ra\pieno)_{opt}$ must pass through a configuration composed by a cluster $E_{B_5}(r^*)$ with the addition of two triangular faces.
	\item[2.] If $\d \in \big(\frac{1}{2},1\big)$, any $\o\in(\vuoto\ra\pieno)_{opt}$ must pass through a configuration composed by a cluster $E_{B_1}(r^*+1)$ with the addition of two triangular faces.
	\ei
\el

\begin{proof}
	We analyze the two cases separately.
	\bi
	\item[1.] Follow the path until it hits $\cV_{A_1^*-3}$. By Lemma \ref{lemma:lower1}, the configuration in this set must be a quasi-regular hexagon with area $6(r^*)^2+10r^*+3$. Since we need not consider any paths that return to the set $\cV_{A_1^*-3}$ afterwards and the path has to cross the set $\cV_{A_1^*-1}$, the path proceeds as follows. Starting from a quasi-regular hexagon with area $A_1^*-3$, a free particle is created giving rise to a configuration with energy $\GK-2\D+2U<\GK$. Before any new particle is created, the energy has to decrease by at least $U$. The unique way to do this is to move the particle towards the cluster and attach it to the quasi-regular hexagon, which lowers the energy to $\GK-2\D+U$. Now it is possible to create another particle at cost $\D$ giving rise to a configuration with energy $\GK-\D+U<\GK$. Again, before creating a new particle, the energy has to decrease by at least $U$. The unique way to do this is to move the particle until it is attached to the cluster, which lowers the energy to $\GK-\D$. Note that this gives us a configuration composed by a cluster $E_{B_5}(r^*)$ with the addition of two triangular faces, as claimed.
	
	\item[2.] We can argue as in the previous case.
	\ei
\end{proof}

\bl{lemma:lower3}
Any $\o\in(\vuoto\ra\pieno)_{opt}$ must reach the energy $\GK$. 
\el

\begin{proof}
By Lemma \ref{lemma:lower2}, we know that any $\o\in(\vuoto\ra\pieno)_{opt}$ must cross a configuration composed by two triangular faces attached to a cluster $E_{B_5}(r^*)$ (resp.\ $E_{B_1}(r^*+1)$) if $\d \in \big(0,\frac{1}{2}\big)$ (resp.\ $\d \in \big(\frac{1}{2},1\big)$). Starting from such a configuration, it is impossible to reduce the energy without lowering the particle number. Indeed, \cite[Theorem 3.22]{AJNT2022} asserts that, for $\d \in \big(0,\frac{1}{2}\big)$ (resp.\ $\d \in \big(\frac{1}{2},1\big)$), the minimal energy in $\cV_{A_1^*-1}$ (resp.\ $\cV_{A_2^*-1}$) is realized (although not uniquely) in such a configuration. Since any further move to increase the particles number involves the creation of a new particle, the energy must reach the value $\GK$.
\end{proof}

\bc{Gamma}
We have
\be{}
\Phi(\vuoto,\pieno)=\GK.
\ee
\ec

\begin{proof}
	Combining Corollary \ref{cor:upperG} and Lemma \ref{lemma:lower3} we obtain the claim.
\end{proof}

\subsection{Structure of the communication level set}

Recalling the two values of the critical area in \eqref{def:criticalarea}, we have the following result.

\bp{prop:ingresso}
The following statements hold:
\bi
\item[1.] Let $\d \in \big(0,\frac{1}{2}\big)$ and $A_1^*=6(r^*)^2+10r^*+6$. Any $\o\in(\vuoto\ra\pieno)_{opt}$ must pass through the set $\cC(A_1^*)=\cK(A_1^*-1)^{fp}$.
\item[2.] Let $\d \in \big(\frac{1}{2},1\big)$ and $A_2^*=6(r^*+1)^2+2(r^*+1)+2$, any $\o\in(\vuoto\ra\pieno)_{opt}$ must pass through the set $\cC(A_2^*)=\cK(A_2^*-1)^{fp}$.
\ei
\ep

\begin{proof}
	We analyze the two cases separately.
	\bi
	\item[1.] By Lemmas \ref{lemma:lower1} and \ref{lemma:lower2}, we can obtain a configuration $\h_0$ with a cluster according to the following cases:
	\begin{enumerate}
		\item[(1)] the two triangular faces form an elementary rhombus which is attached to one of the longest sides of the quasi-regular hexagonal cluster, namely the resulting configuration is in $\tilde\cS(A_1^*-1)$ (see Figure \ref{addtriangularunitK});
		\item[(2)] the two triangular faces are attached to one of the longest sides of the quasi-regular hexagonal cluster at triangular lattice distance 2, namely the resulting configuration is in $\tilde\cD(A_1^*-1)$ (see Figure \ref{addtriangularunitK});
		\item[(3)] the two triangular faces are attached to the same side of the quasi-regular hexagonal cluster at triangular lattice distance greater than 2 (see Figure \ref{addtriangularunitK});
		\item[(4)] the two triangular faces are attached to two different sides of the quasi-regular hexagonal cluster (see Figure \ref{addtriangularunitK});
		\item[(5)] the two triangular faces form an elementary rhombus which is attached to one of the sides, other than the longest, of the quasi-regular hexagonal cluster;
		\item[(6)] the two triangular faces are attached at triangular lattice distance 2 to the same side, other than the longest, of the quasi-regular hexagonal cluster;
        \item[(7)] the two triangular faces form an elementary rhombus which is attached to one of the sides, but the direction of the elementary rhombus is towards the outer direction of the cluster.
	\end{enumerate}
	
	\begin{figure}
		\centering \includegraphics[scale=0.4]{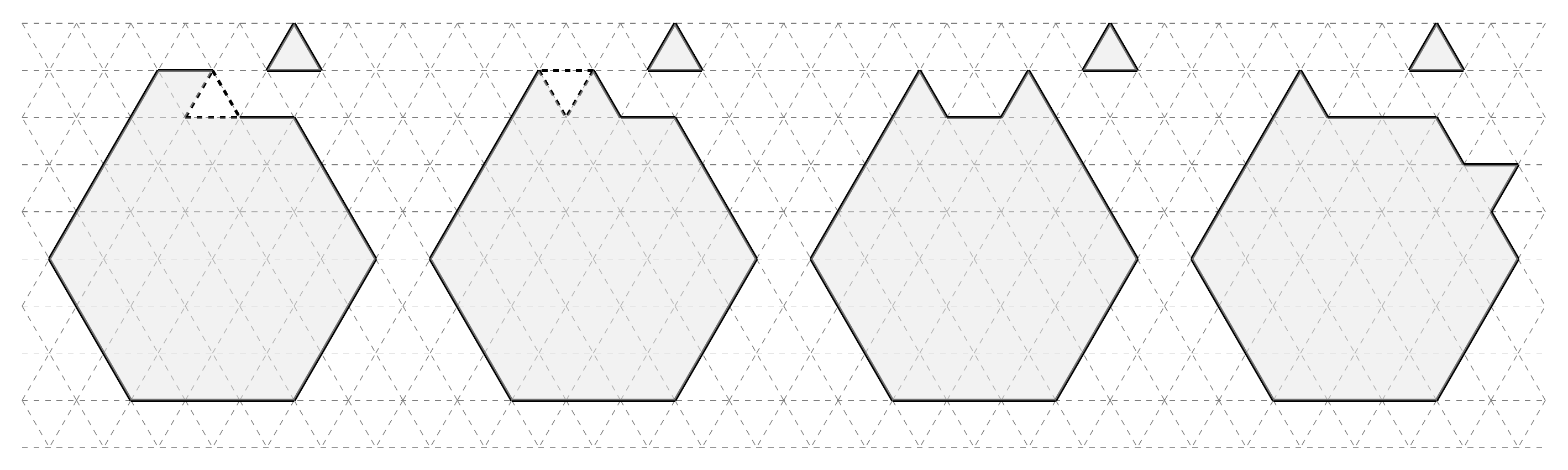} 
		\caption{From the left to the right we depict in light grey the clusters described in cases (1)--(2)--(3)--(4). In dashed dark grey, we depict the future position of the free particle to cover the angle of $\frac{5}{3} \pi$ in the first two cases.} \label{addtriangularunitK}
	\end{figure}

	Note that in all these cases the cluster has minimal perimeter, indeed it has the same perimeter as a standard hexagon with the same area. Moreover, in all these cases the configuration $\h_0$ has energy $\GK-\D$. We will prove that every $\o\in{(\vuoto\ra\pieno)_{opt}}$ crosses a configuration in $\cC(A_1^*)$. Since we need not consider any paths that return to the set $\cV_{A_1^*-2}$ afterwards and the energy can increase by at most $\Delta$ in order to have an optimal path, there are only the following possibilities:
	\begin{itemize}
		\item[A.] a free particle enters $\L$;
		\item[B.] a particle is detached from the cluster; 
		\item[C.] a particle moves at cost $U$ without detaching from the cluster. 
	\end{itemize}
	
	{\bf Case A.} We may assume that the free particle does not exit from $\L$, otherwise we can iterate this argument for a finite number of steps since the path has to reach $\pieno$. Let $\h_1=\h_0^{fp}$. Since $H(\h_1)=\G^{\text{K-Hex}}$, in order to have an optimal path the energy cannot increase. Thus the unique admissible moves are the movement of the free particle at zero cost and the attachment of the particle to the cluster. We may assume that the particle attaches to the cluster, otherwise we can iterate this argument. 
	
	In cases (1) and (2), note that $\h_1$ contains an internal angle of $\frac{5}{3}\pi$, thus  we consider the configuration $\h_2$ obtained from $\h_1$ by attaching the free particle to cover the internal angle of $\frac{5}{3}\pi$ of the cluster (see Figure \ref{addtriangularunitK}).  
	Thus the energy decreases by $2U$ and therefore it is possible to create a new particle without exceeding the energy value $\GK$. Indeed, let $\h_3$ be the configuration obtained from $\h_2$ by creating a new particle, thus we obtain:
	\be{sottoG}
	\ba{ll}
	H(\h_3)&=(H(\h_3)-H(\h_2))+(H(\h_2)-H(\h_1))+H(\h_1) \\
	&=\G^{\text{K-Hex}}+\D-2U<\G^{\text{K-Hex}}.
	\ea
	\ee
	From now on the path proceeds as the reference path $\o^*$ without exceeding the energy value $\GK$. Note that the path crosses the set $\cC(A_1^*)$ in the configuration $\h_1$.
	
	In cases (5) and (6), since $\h_1$ contains an internal angle of $\frac{5}{3}\pi$, it is possible that the free particle attaches to the cluster at cost $-2U$. If this occurs, we can derive \eqref{sottoG} as before, but we show that now it is not possible to reach $\pieno$ without exceeding $\GK$ unless the path reaches a configuration $\h_0$ as in cases (1) or (2). Any side $s$ of $E_{B_5}(r^*)$, other than the longest, has length $r^*+1$, so the bar $B$ with the larger base $r^*+1$ has cardinality $l=2r^*-1$. We can write
	\be{}
	\o=(\vuoto,...,\tilde\h,...,\h_0,\h_{1},...,\h_{i_l},...,\h_{i_{l+1}},...,\h_{i_{l+2}},\bar\h,...,\pieno),
	\ee
	where $\tilde\h=\cE_{B_5}(r^*)$, $\h_0$ and $\h_{1}$ are as above, $\h_{i_l}$ is the configuration obtained after filling the new bar $B$ and creating a free particle, $\h_{i_{l+1}}$ is the configuration obtained from $\h_{i_l}$ by attaching the free particle and afterwards creating another free particle, and $\h_{i_{l+2}}$ is the configuration obtained from $\h_{i_{l+1}}$ by attaching the free particle to the cluster. Finally, let $\bar\h$ the configuration obtained from $\h_{i_{l+2}}$ by creating a free particle. Thus we obtain the following contradiction:
	\be{eccedo}
	\ba{ll}
	H(\bar\h)&=(H(\bar\h)-H(\h_{i_{l+2}}))+(H(\h_{i_{l+2}})-H(\h_{i_{l+1}}))+(H(\h_{i_{l+1}})-H(\h_{i_l})) \\
	&\quad+(H(\h_{i_l})-H(\h_{1}))+H(\h_1)\\
	&=\D-U+(\D-U)+((2r^*-3)\Delta+(4-3r^*)U)+\G^{\text{K-Hex}}>\G^{\text{K-Hex}}.
	\ea
	\ee
	Therefore, starting from the configuration $\h_{1}$, after attaching the protuberance at cost $-2U$ the path cannot sequentially create and attach a particle to the cluster: this follows from \eqref{eccedo}. Thus the path has to further lower the energy before reaching the configuration $\h_{i_{l+2}}$. If the path reaches a configuration $\xi$ such that $n(\xi)=0$, i.e., $\xi$ has no free particle, a free particle has been attached and the energy lowered by $2U$ at most, but this does not suffice due to \eqref{eccedo}. But there are no moves that further lower the energy. If the path reaches a configuration $\xi$ such that $n(\xi)=1$, then the unique way to lower the energy is to attach the free particle at cost $-2U$ or $-U$, but again this does not suffice due to \eqref{eccedo}. Since the path $\o$ has to reach $\pieno$ and therefore the number of particles has to increase, the unique possibility in order to have an optimal path is that the path $\o$ comes back to the configuration $\h_{0}$. Thus we are done as claimed before.

	\begin{figure}
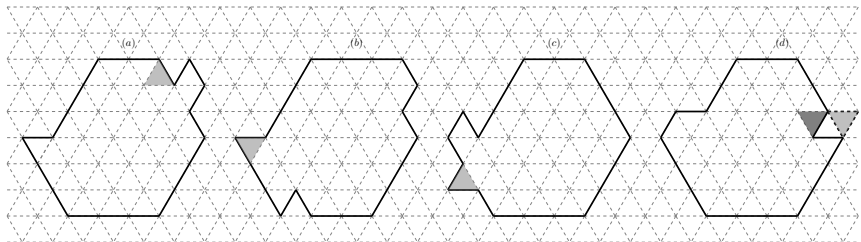

		\centering 
		\includestandalone[scale=0.4, mode=image|tex]{trenino1}
		\caption{In (a) we depict an example of configuration with $\h_0$ as in case (5) and a protuberance attached at cost $-U$, in which we highlight in grey the triangular face that moves towards the elementary rhombus at cost $U$, while in (b) we depict the final configuration of these moves, in which we highlight in grey the protuberance that has to be detached. In (c) and (d) we depict an example of configuration with $\h_0$ as in case (3) and a protuberance attached at cost $-U$: in (c) we highlight in grey the single protuberance, while in (d) we highlight in light grey the free particle and in dark grey the last triangular face that has been moved towards left.}
		\label{trenino1}
	\end{figure}

	If the particle attaches at cost $-U$, motions of particles at cost $U$ can take place. The unique possibility to first move a particle at cost $U$ is to attach it to the elementary rhombus (see Figure \ref{trenino1}(a)), otherwise it is possible to move the last attached protuberance, but in this case we can iterate the argument. All the configurations that are crossed during these motions have energy either $\GK$ or $\GK-U$. Since the energy has to further lower in order to create a new particle and reach $\pieno$, the unique possibility is to detach a protuberance at cost $U$ and attach it to cover the unique internal angle of $\frac{5}{3}\pi$ (see Figure \ref{trenino1}(b)). Thus we can argue as before for cases (5) and (6).
	
	In cases (3), (4) and (7) note that the configuration $\h_1$ does not contain an internal angle of $\frac{5}{3}\pi$, thus the free particle can attach only at cost $-U$. Let $\h_2$ be the configuration obtained from $\h_1$ by attaching the free particle to the cluster at cost $-U$, thus $H(\h_2)=\GK-U$. Note that now the unique admissible moves are those at cost $U$ at most, thus it is not possible to create a new particle before further lowering the energy.
	
	 In cases (3) and (7), since the path has to reach $\pieno$, the unique possibility is to move a protuberance $T$ in such a way it forms an angle of $\frac{5}{3}\pi$ with another protuberance.
 If the two protuberances in configuration $\h_0$ are not on the longest side, we are left to analyze case A for a configuration $\h_2$ as in cases (5) or (6) and therefore we can argue as before. If the two protuberances in configuration $\h_0$ are on the longest side, we obtain a configuration composed by a cluster as in case (1) or (2) with the addition of a protuberance. As explained before, it is possible that motions at cost $U$ take place. All the configurations that are crossed during these motions have energy either $\GK$ or $\GK-U$. At the end of these motions, there are only the two following possibilities: there is a unique cluster with an internal angle of $\frac{5}{3}\pi$ and either a free particle (see Figure \ref{trenino1}(d)) or a single protuberance (see Figure \ref{trenino1}(c)). In the first case, the configuration that is obtained is in $\cC(A_1^*)$. In the latter case, since the energy has to further lower in order to create a new particle and reach $\pieno$, the unique possibility is to detach the single protuberance at cost $U$ and attach to the cluster at cost $-2U$. When the protuberance is detached the path crosses the set $\cC(A_1^*)$.
	
	In case (4), if the third protuberance is attached in such a way all the three protuberances are attached to different sides, then the unique admissible moves are detaching a protuberance. Thus we can iterate this argument for a finite number of steps, since the path has to reach $\pieno$. We are left to consider the case in which at least two protuberances are attached to the same side. We can argue as above.

	{\bf Case B.} Let $\h_1$ be the configuration obtained from $\h_0$ by detaching a particle from a cluster. Since $H(\h_0)=\GK-\D$ and the path $\o$ has to be optimal, the energy can increase by $U$ at most. Thus, only a protuberance can be detached. After that, only moves with cost 0 at most are admissible. Since the path has to reach $\pieno$ and therefore the free particle cannot move for infinite time at zero cost, the unique possibility is to attach the free particle to the cluster. Thus we obtain a configuration that is analogue to $\h_0$ and we can iterate this argument for a finite number of steps, until we come back to case A.
	
	{\bf Case C.} Note that this case is admissible only for configurations $\h_0$ as in cases (1), (2), (5) or (6). Since $H(\h_0)=\GK-\D$ and the path $\o$ has to be optimal, the energy can increase by $U$ at most.
	
	In cases (1) or (2), starting from $\h_0$, all the configurations that can be obtained without exceeding the energy value $\GK$ are in the set $\cK(A_1^*-1)$: this directly follows from the definition of that set given in \eqref{proto} since no particle is detached from the cluster. From now on, since the energy of the last configuration is $\GK-\D$, it is possible either to create a free particle, or to detach a protuberance at cost $U$, or to move a particle at cost $U$. In the first case, we can conclude as in case A. Note that the path $\o$ crosses the set $\cC(A_1^*)$ when the free particle is created. In the second case, we can conclude as in case B. In the latter case, we can iterate this argument for a finite number of steps.

	In cases (5) or (6), we can argue as for the cases (1) or (2). Indeed, the same kind of motions can take place, but when the free particle is created, we can argue as in the case A for $\h_1$ as in cases (5) or (6). This concludes the proof.

	\item[2.] The proof of this case is similar to the previous one.
	\ei
\end{proof}

\br{remark:ingresso}
We want to emphasize that this result is different from \cite[Proposition 2.3.7]{BHN}. Indeed, on the square lattice the authors were able to prove that any optimal path from the metastable to the stable state reaches a square or quasi-square shape, then a protuberance is attached and finally a free particle enters the box. However, Lemma \ref{lemma:lower2} and Proposition \ref{prop:ingresso} do not suffice to characterize the entrance in the gate. Indeed, several mechanism to enter the gate appear on the hexagonal lattice. Clearly, one of these possibilities is to add a free particle starting from a configuration in $\cK(A_i^*-1)$, with $A_i^*\in\{A_1^*,A_2^*\}$. But there are many other ways to enter $\cK(A_i^*-1)^{fp}$. For example, suppose that $0<\d<\frac{1}{2}$ and an optimal path $\o:\vuoto\ra\pieno$ crosses a configuration $\h$ of the type described in case (3) in the proof of Proposition \ref{prop:ingresso} with the addition of a free particle. Starting from $\h$, it is possible that the free particle is attached to the cluster in such a way that it forms an elementary rhombus together with a triangular face already attached. Thus the energy reaches the value $\GK-U$. Thus, it is possible to move the other triangular face at cost $U$ and, when it is detached, the path $\o$ crosses a configuration in $\tilde \cS(A_1^*-1)^{fp}\subset\cK(A_1^*-1)^{fp}$, but the path does not cross the set $\cK(A_1^*-1)$. With this example we want to put the attention on the fact that several mechanisms to enter the gate appear due to the particular shape of the lattice. Indeed, on the square lattice it does not matter which side the protuberance is attached to because it is possible to move it along the side at zero cost.
\er

\section{Recurrence property}\label{recurrencepropertyproof}
The goal of this Section is to prove Theorem \ref{prop:recurrence_property}. Recall \eqref{stab} for the definition of stability level.
The following theorem states that every configuration of $\mathcal{X}$ different from $\vuoto$ and $\pieno$ has a stability level $\Delta+U$ at most. 
\bp{prop:stability_lower}
	Let $\eta \in \mathcal{X}$ be a configuration such that $\eta \not \in \{\vuoto,\pieno\}$, then $V_{\eta} \leq \Delta+U$.
\ep
An immediate consequence of Proposition \ref{prop:stability_lower} is that the only configurations with a stability level greater than $\Delta+U$ are $\vuoto$ and $\pieno$, as reported in Theorem \ref{prop:recurrence_property}. 
The proof of Proposition \ref{prop:stability_lower} is divided in two steps. First of all, in Section \ref{section5.1} we prove that the configurations with peculiar geometrical properties has a stability level smaller than or equal to $\D+U$ (see Lemmas \ref{nofreeparticles}-\ref{nointeracting}), and then in Section \ref{section5.2} we show that all configurations, different from $\vuoto$ and $\pieno$, has a stability level smaller than or equal to $\Delta+U$, i.e., $\mathcal{X}_{\D+U} \setminus \{\vuoto, \pieno \}=\emptyset$ (see Lemma \ref{zry}). 

\subsection{Configurations with stability level \texorpdfstring{$\D+U$}{TEXT} at most}\label{section5.1}

Recall \eqref{xv} for the definition of $V$-irreducible states. 
In this Section, we emphasize the different stability level for configurations depending on their particular geometrical properties. For the proof of the lemmas we refer to Section \ref{dimlemmi}. The following Lemma characterizes the configurations in $\cX_0$.

\bl{nofreeparticles}
Any configuration $\h\in\cX_{0}$ has no free particles.
\el

In order to state the following Lemmas, we need the following definition. Moreover, recall Definition \ref{def:hole_in_polyiamond} for the definition of a hole.

\bd{}
Two clusters are called \emph{interacting} if their lattice distance is two. Otherwise, two clusters are called \emph{non-interacting} if its lattice distance is strictly greater than 2. 
\ed

\bl{stablevU}
	If a configuration $\sigma$ contains a cluster with an internal angle of $\frac{1}{3}\pi$ and no free particles, no holes and no interacting clusters, then it has a stability level smaller than or equal to $U$, i.e., $\sigma \not \in \mathcal{X}_{U}$.
\el

\bl{stablevdelta}
	If a configuration $\sigma$ contains a cluster with an internal angle of $\frac{5}{3}\pi$ and no free particles, no holes and no interacting clusters, then it has a stability level smaller than or equal to $\Delta$, i.e., $\sigma \not \in \mathcal{X}_\Delta$.
\el

\bl{stablev2}
	If a configuration $\sigma$ contains a cluster with an internal angle of $\frac{4}{3}\pi$ and no free particles, no holes and no interacting clusters, then it has a stability level smaller than or equal to $2\Delta-U$, i.e., $\sigma \not \in \mathcal{X}_{2\Delta-U}$.
\el

The next lemma investigates the case in which a configuration contains two interacting clusters or a cluster with a hole.

\bl{nointeracting}
	If a configuration $\sigma$ contains two interacting clusters or a cluster with a hole, then it has a stability level smaller than or equal to $\Delta+U$, i.e., $\sigma \not \in \mathcal{X}_{\Delta+U}$.
\el

\subsection{Identification of configurations in \texorpdfstring{$\mathcal{X}_{\Delta+U}$}{TEXT}}\label{section5.2}
In Section \ref{section5.1}, we established that the configurations with particular geometrical properties has a stability level $\D+U$ at most. The configurations, that do not satisfy Lemma \ref{nofreeparticles} and Lemma \ref{nointeracting}, has no free particle and no interacting clusters. Moreover, the configurations, that do not satisfy Lemmas \ref{stablevU}-\ref{stablev2}, contain clusters with internal angles of $\pi$ and $\frac{2}{3}\pi$ only. Thus, the clusters contained in these configurations have an hexagonal shape. Now, we partition the set of remaining configurations, different from $\vuoto$ and $\pieno$, into three subsets $Z, R, Y$ and we prove that also these configurations has a stability level smaller than or equal to $\D+U$. Thus, it follows that if there exists a configuration with a stability level strictly greater than $\D+U$, then it is $\vuoto$ or $\pieno$.

\begin{figure}
	\centering
	\includegraphics[scale=0.7]{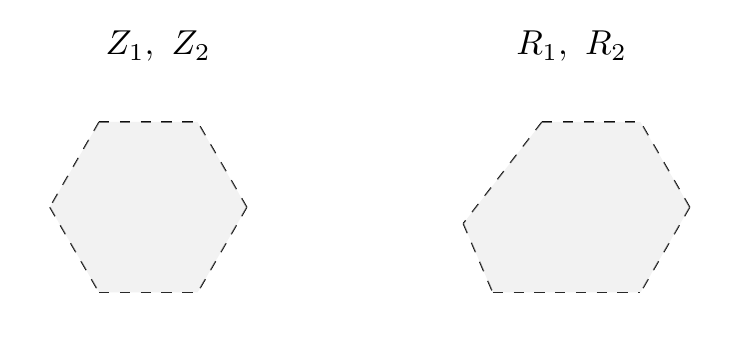}
	\caption{On the left-hand side we depict an example of the cluster in a configuration of $Z$, while on the right-hand side an example of the cluster in a configuration of $R$.}
	\label{configurationZR}
\end{figure}

$Z$ is the set of configurations consisting of a single quasi-regular hexagonal cluster (see Figure~\ref{configurationZR} on the left-hand side). More precisely, $Z=Z_1 \cup Z_2$, where:
\begin{itemize}
	\item $Z_1$ is the collection of configurations such that there exists only one cluster with shape $E_{B_m}(r) \subset \Lambda$ with $r \leq r^*$ and $m \in \{0,1,2,3,4,5\}$;
	\item $Z_2$ is the collection of configurations such that there exists only one cluster with shape $E_{B_m}(r) \subset \Lambda$ with $r \geq r^*+1$ and $m \in \{0,1,2,3,4,5\}$.
\end{itemize}

We define the set $R$ to be the set of configurations consisting of a single hexagonal cluster (see Figure~\ref{configurationZR} on the right-hand side). Formally, $R=R_1 \cup R_2$, where:
\begin{itemize}
	\item $R_1$ is the collection of configurations such that there exists only one cluster with hexagonal shape $E \subset \Lambda$ such that it contains the greatest quasi-regular hexagon with radius $r \leq r^*$; 
	\item $R_2$ is the collection of configurations such that there exists only one cluster with hexagonal shape $E \subset \Lambda$ such that it contains the greatest quasi-regular hexagon with radius $r \geq r^*+1$.
\end{itemize}

\begin{figure}
	\centering
	\includegraphics[scale=0.7]{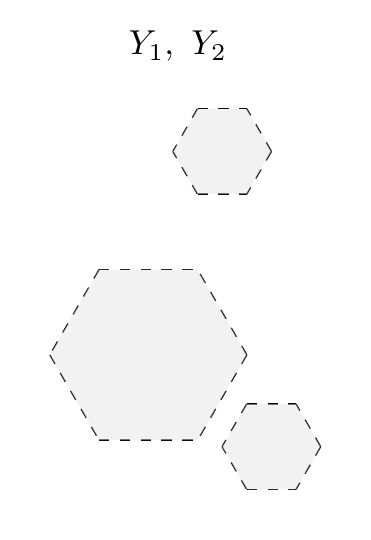}
	\caption{An example of clusters in a configuration of $Y$.}
	\label{configurationY}
\end{figure}

The set $Y$ contains all configurations with more than one 
hexagonal cluster of the types in $Z_1$, $Z_2$, $R_1$, $R_2$ (see Figure~\ref{configurationY}). More precisely, we have $Y=Y_1 \cup Y_2$, where: 
\begin{itemize}
	\item $Y_1$ is the collection of configurations such that there exists a family of non-interacting clusters with hexagonal shape such that it contains the greatest quasi-regular hexagon with radius $r \leq r^*$; 
	\item $Y_2$ is the collection of configurations such that there exists a family of clusters with at least one having hexagonal shape containing the greatest quasi-regular hexagon with radius $r \geq r^*+1$.
\end{itemize}
In other words $Y_1$ contains a collection of clusters of the same type of those in $Z_1$ or $R_1$, and $Y_2$ contains a collection of clusters where at least one is of the same type of those in $Z_2$ or $R_2$.

\bl{zry}
	If $\s\in Z\cup R\cup Y$, then $V_{\sigma}\leq \D+U$.
\el

We refer to Section \ref{dimlemmi} for the proof of this lemma.

\subsection{Proof of Lemmas}\label{dimlemmi}

\begin{proof*}{\it of Lemma~\ref{nofreeparticles}}
	If $\h$ has a free particle, then $\h$ is obviously $0$-reducible, i.e., its stability level is $0$ and therefore $\h\notin \cX_{0}$. Indeed, the reducing path is immediately obtained by bringing the free particle outside $\L$ or attaching it to a cluster.
\end{proof*}

\begin{figure}
	\centering \includegraphics[scale=0.6]{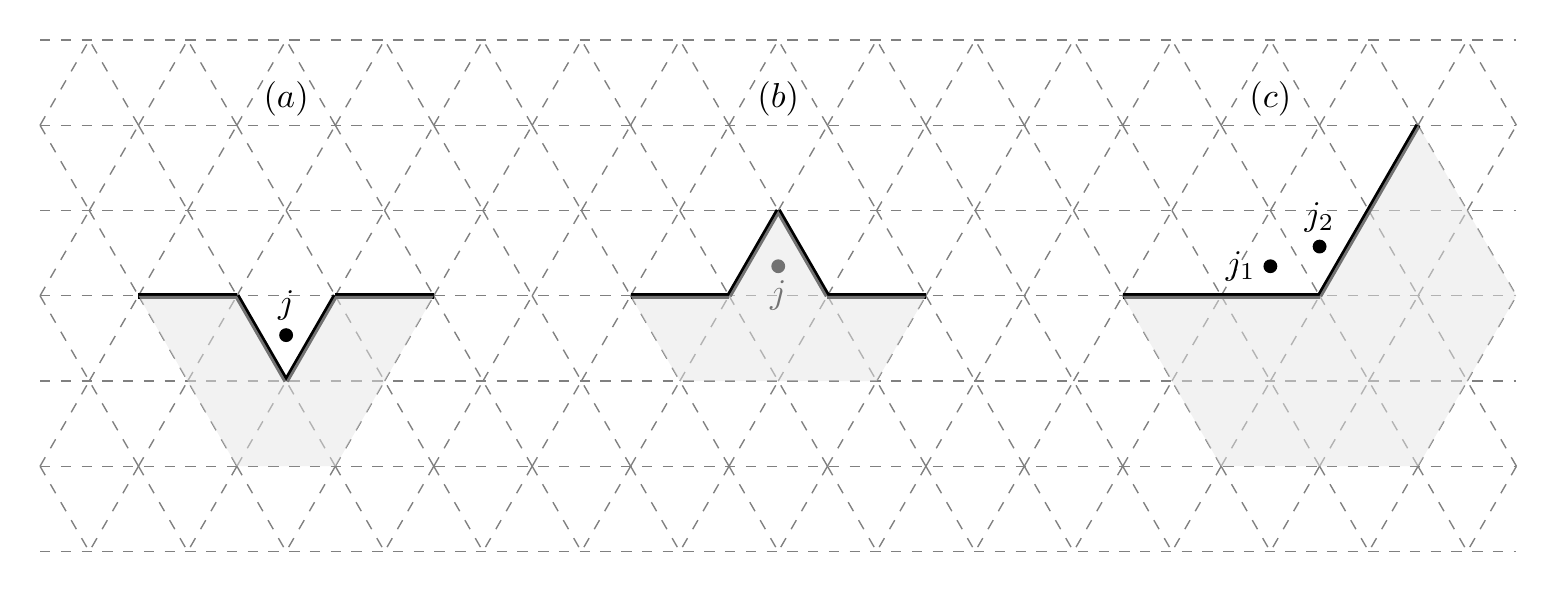} 
	\caption{On the left-hand side (resp.\ center) we depict the site $j$ when $\sigma$ has an internal angle of $\frac{5}{3}\pi$ (resp.\ $\frac{1}{3}\pi$). On the right-hand side we depict the two sites $j_1,j_2$ when $\sigma$ has an internal angle of $\frac{4}{3}\pi$.} \label{anglesconf}
\end{figure}

\begin{proof*}{\it of Lemma~\ref{stablevU}}
	Let $\s$ be a configuration as in the statement and let $C(\sigma)$ be a cluster with an internal angle $\alpha=\frac{1}{3}\pi$. Let $j$ be a site such that $\sigma(j)=1$ and that belongs to the closed triangular face of $C(\sigma)$ intersecting its boundary in two edges (see Figure~\ref{anglesconf}(b)). We define $\eta$ as the configuration obtained from $\sigma$ by detaching the particle in $j$ and then moving it outside $\L$. Note that it is possible to bring the particle outside $\L$ since $\sigma$ does not contain clusters with holes or interacting clusters. We construct a path $\omega:\s\ra\h$ as 
	\begin{equation}
		\omega=(\sigma, \xi_1, \xi_2, ..., \xi_n, \eta),
	\end{equation}
	where $\x_1$ is the configuration obtained from $\s$ by moving the particle at site $j$ in one of the two empty nearest-neighbor sites. The cost of this move is $U$. Then the particle, after being detached, can be brought outside $\L$ passing through the configurations $\x_2,...,\x_n$, possibly interacting with other clusters. We need to bring the particle outside $\L$ because the energy does not necessarily decreases by $2U$ when the particle interacts with the clusters during the motion. Note that the energy of the configurations $\x_2,...,\x_n$ is $H(\xi_1)$ at most. Thus we obtain 
	\begin{equation}
		H(\eta)-H(\sigma)=U-\Delta<0,
	\end{equation}
	where the inequality follows from the condition $\Delta>U$. Thus, $\eta$ belongs to $\mathcal{I}_{\sigma}$ and $V_{\sigma} \leq U$.
\end{proof*}

\begin{proof*}{\it  of Lemma~\ref{stablevdelta}}
	Let $\s$ be a configuration as in the statement and let $C(\sigma)$ be a cluster with an internal angle $\alpha=\frac{5}{3}\pi$. Let $j$ be the site at distance one to a site in $C(\sigma)$ such that $\sigma(j)=0$ and that belongs to the closed triangular face intersecting the boundary of $C(\sigma)$ in two or more edges (see Figure \ref{anglesconf}(a)). We define $\eta$ as the configuration obtained by $\sigma$ after creating a particle and then attaching it in the site $j$. Note that it is possible to bring the particle from the boundary of $\L$ towards the site $j$ since $\sigma$ does not contain clusters with holes or interacting clusters. We construct a path $\omega$ connecting $\sigma$ and $\eta$ as 
	\begin{equation}
		\omega=(\sigma, \xi_1, \xi_2, ..., \xi_n, \eta),
	\end{equation}
	where $\x_1$ is the configuration obtained from $\sigma$ by creating a particle in $\partial^-\L$ at cost $\D$. Then, this particle moves towards the cluster $C(\sigma)$, passing through the configurations $\x_2,...,\x_n$, until it is attached in the site $j$ at cost $-2U$ giving rise to the configuration $\h$. Note that the energy of the configurations $\x_2,...,\x_n$ is $H(\xi_1)$ at most. Thus we obtain
	\begin{equation}
		H(\eta)-H(\sigma)=\Delta-2U<0,
	\end{equation}
	where the inequality follows from $\Delta<\frac{3}{2}U$. Thus, $\eta$ belongs to $\mathcal{I}_{\sigma}$ and $V_{\sigma} \leq \Delta$.
\end{proof*}

\begin{proof*}{\it of Lemma~\ref{stablev2}}
	Let $\s$ be a configuration as in the statement and let $C(\sigma)$ be a cluster with an internal angle $\alpha=\frac{4}{3}\pi$. Let $j_1$, $j_2$ be two sites such that $\sigma(j_1)=\sigma(j_2)=0$, $d(j_1,j_2)=1$ and let each of them belong to one closed triangular face intersecting the boundary of $C(\sigma)$ in one edge (see Figure \ref{anglesconf}(c)). 
	We define $\eta$ as the configuration obtained by $\sigma$ after the following sequence of moves: creation of a particle and movement of it until it is attached in the site $j_1$; creation of another particle and movement of it until it is attached in the site $j_2$. Note that it is possible to bring particles from the boundary of $\L$ towards the sites $j_1$ and $j_2$ since $\sigma$ does not contain clusters with holes. We construct a path $\omega$ connecting $\sigma$ and $\eta$ as 
	\begin{equation}
		\omega=(\sigma, \xi_{i_1}, \xi_{i_2}, ..., \xi_{i_n},\x, \x_{j_1},...,\x_{j_m}, \eta),
	\end{equation}
	where $\x_{i_1}$ is the configuration obtained from $\sigma$ by creating a particle in $\partial^-\L$ at cost $\D$. Then, this particle moves towards the cluster $C(\sigma)$, passing through the configurations $\x_{i_2},...,\x_{i_n}$, until it is attached in the site $j_1$ at cost $-U$ giving rise to the configuration $\x$. Note that the energy of the configurations $\x_{i_2},...,\x_{i_n}$ is $H(\xi_{i_1})$ at most. The configuration $\x_{j_1}$ is obtained from $\x$ by creating a particle in $\partial^-\L$ at cost $\D$. Then, this particle moves towards the cluster $C(\sigma)$, passing through the configurations $\x_{j_2},...,\x_{j_m}$, until it is attached in the site $j_2$ at cost $-2U$ giving rise to the configuration $\h$. Note that the energy of the configurations $\x_{j_2},...,\x_{j_m}$ is $H(\xi_{j_1})$ at most. Thus we obtain
	\begin{equation}
		H(\eta)-H(\sigma)=2\Delta-3U<0,
	\end{equation}
	where the inequality follows from $\Delta <\frac{3}{2} U$. Thus, $\eta$ belongs to $\mathcal{I}_{\sigma}$ and $V_{\sigma} \leq 2\Delta-U$.
\end{proof*}

\begin{proof*}{\it of Lemma~\ref{nointeracting}}
	We analyze the configuration $\sigma$ starting from the clusters with minimal distance to the boundary of $\Lambda$. If the first clusters $C_1(\sigma)$ and $C_2(\sigma)$, according to this minimal distance, are interacting, we consider the shared vertex $v:=C_1(\sigma) \cap C_2(\sigma)$ on the triangular lattice, see the first and the second pictures in Figure \ref{interacting}. We let into $\Lambda$ a particle and we call $\s_1$ this new configuration. Then this new particle moves until it is attached to $v$ giving rise to the configuration $\s_2$. In this way, there are two possibilities:
	\begin{itemize}
		\item[(i)] the triangular face of this particle shares an edge with $C_1(\sigma)$, an edge with $C_2(\sigma)$ and contains $v$, see the second picture in Figure \ref{interacting};
		\item[(ii)] the triangular face of this particle contains $v$ and shares an edge either with $C_1(\sigma)$ or $C_2(\sigma)$, see the first picture in Figure \ref{interacting}.
	\end{itemize}

	\begin{figure}
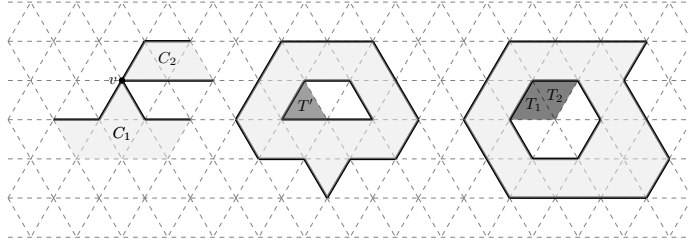

		\centering 
		\includestandalone[scale=0.6, mode=image|tex]{interactingclusters}
		\caption{Here we represent clusters in light grey. In the first two pictures we depict two interacting clusters $C_1$ and $C_2$, while in the third one (resp.\ fourth one) we depict an example of configuration belonging to case 2 (resp.\ case 3) in the proof of Lemma \ref{nointeracting}. }
		\label{interacting}
	\end{figure}
	
	{\bf Case (i).} We have
	\begin{equation}
		\begin{array}{ll}
			H(\sigma_1)-H(\sigma)&= \Delta, \\
			H(\sigma_2)-H(\sigma_1)&= -2U
		\end{array}
	\end{equation}
	and therefore
	\begin{equation}
		H(\sigma_2)-H(\sigma)=[H(\sigma_2)-H(\sigma_1)]+[H(\sigma_1)-H(\sigma)]=-2U+\Delta<0.
	\end{equation}
	Thus, the stability level of $\sigma$ in this case is $V_{\sigma}=\Delta$.
	
	{\bf Case (ii).} We may assume without loss of generality that the triangular face $T$ of the new particle shares an edge with $C_1(\sigma)$. We have that 
	\begin{equation}
		\begin{array}{ll}
			H(\sigma_1)-H(\sigma)&= \Delta, \\
			H(\sigma_2)-H(\sigma_1)&= -U
		\end{array}
	\end{equation}
	and therefore 
	\begin{align}
		H(\sigma_2)-H(\sigma)=-U+\D>0.
	\end{align}
	since $\D>U$. Thus, the energy has to further lower. 
	We define the configuration $\s_3$ as the one obtained from $\s_2$ by creating a new particle. Then this particle moves towards the cluster until it is attached close to the triangular face $T$ in such a way it is attached also to the cluster $C_2(\s)$. 
	This configuration is called $\s_4$. Thus, we have
	\begin{equation}
		\begin{array}{ll}
			H(\sigma_3)-H(\sigma_2)&= \D, \\
			H(\s_4)-H(\s_3)&=-2U. 
		\end{array}
	\end{equation}
	Follows that
	\begin{equation}
		\begin{array}{ll}
			H(\sigma_4)-H(\sigma)&=[H(\sigma_4)-H(\sigma_3)]+[H(\sigma_3)-H(\sigma_2)]\\
			&\quad+[H(\sigma_2)-H(\sigma_1)]+[H(\sigma_1)-H(\sigma)] \\
            &=-3U+2\Delta<0
		\end{array}
	\end{equation}
	and the stability level of $\sigma$ in this case is $V_{\sigma}=2\Delta-U$.
	
	Thus we conclude that the stability level for a configuration with two interacting clusters is $\max\{\D,2\D-U\}=2\D-U$.
	
	Next, suppose that the first clusters are not interacting. Let $C(\sigma)$ be the first cluster of $\sigma$ with a hole.
	Consider one of the empty triangular faces in the hole that share at least an edge with the cluster. There are three cases:
	\begin{itemize}
		\item[1.] the empty triangular face shares three edges with the cluster;
		\item[2.] the empty triangular face shares two edges with the cluster, which we represent with the dark grey triangular face $T'$ in the third picture in Figure \ref{interacting};
		\item[3.] the empty triangular face shares only one edge with the cluster, which we represent with the dark grey triangular face $T_1$ in the fourth picture in Figure \ref{interacting}.
	\end{itemize}
	In the first case, we move the empty triangular face until it reaches the internal boundary of the cluster. Since every triangular face in the internal boundary of the cluster shares at least an edge with the empty triangular faces outside the cluster, then
	\begin{align}
		H(\eta)-H(\sigma) \leq -U,
	\end{align}
	where $\eta$ is the configuration obtained from $\sigma$ by exchanging the empty triangular face of the hole with a triangular face on the internal boundary of the cluster. Thus $V_\sigma=0$.
	
	In the second case, as before, we move the empty triangular face $T'$ until it reaches the internal boundary of the cluster giving rise to the configuration $\s_1$. If $\s$ and $\s_1$ can be connected via one step of the dynamics, then the energy value remains the same. Otherwise, during the first step, the energy increases by $U$, indeed the empty triangular face $T'$ can be detached from the other empty triangular face by breaking two bonds and creating only a new bond (see the third picture in Figure \ref{interacting}). Thus in both cases it holds that $H(\s_1)-H(\s)\leq U$. Moreover, every triangular face in the internal boundary of the cluster shares one or two edges with the empty triangular face outside the cluster. 
	\begin{itemize}
		\item If there exists a triangular face $T$ in the internal boundary of the cluster $C(\sigma)$ with two shared edges with some empty triangular faces, then, denoting by $\eta$ the configuration obtained from $\sigma_1$ by exchanging the empty triangular face of the hole with $T$, we have
		\begin{align}
			& H(\eta)-H(\sigma_1)= -2U, \\
			& H(\eta)-H(\sigma)=[H(\eta)-H(\sigma_1)]+[H(\sigma_1)-H(\sigma)]\leq-U.
		\end{align}
		Thus $V_\sigma=U$.
		
		\item Otherwise, if each triangular face in the internal boundary has only one shared edge with an empty triangular face outside cluster, then we have 
		\begin{equation}
			H(\sigma_2)-H(\sigma_1)=-U,
		\end{equation} 
		where $\sigma_2$ is the configuration obtained from $\sigma_1$ by exchanging the empty triangular face of the hole with a triangular face in the internal boundary of the cluster. Thus, we obtain $H(\sigma_2)=H(\sigma)$, and by construction $\sigma_2$ has an internal angle of $\frac{5}{3}\pi$ (see Figure \ref{anglesconf}(a)). We define a configuration $\sigma_3$ obtained from $\sigma_2$ by getting in $\Lambda$ a new particle, and we define $\sigma_4$ from $\sigma_3$ by moving and attaching this particle to cover the internal angle of $\frac{5}{3}\pi$. We have
		\begin{equation}
			\begin{array}{ll}
				H(\sigma_3)-H(\sigma_2) &= \Delta, \\
				H(\sigma_4)-H(\sigma_3) &=-2U
			\end{array}
		\end{equation}
		and therefore
		\begin{equation}
       \begin{array}{ll}
			H(\sigma_4)-H(\sigma)&=[H(\sigma_4)-H(\sigma_3)]+[H(\sigma_3)-H(\sigma_2)]+[H(\sigma_2)-H(\sigma)] \\
            &\leq-2U+\Delta<0. 
            \end{array}
		\end{equation}
		Thus $V_\sigma=\Delta$.
	\end{itemize}
	
	In the third case, the empty triangular face $T_1$ has only one shared edge with the cluster, so there exists another empty triangular face $T_2$ in the hole that is connected with $T_1$. We move the two empty triangular faces until they reach the internal boundary of the cluster giving rise to the configuration $\s_1$. If $\s$ and $\s_1$ can be connected via one step of the dynamics, then the energy value increases by $U$. Otherwise, during the first step, the energy increases by $2U$, indeed the empty triangular face moves $T_1$ from the other empty triangular face by breaking two bonds (see the fourth picture in Figure \ref{interacting}). Thus in both cases we have that 
	\begin{equation}\label{stab2U_nobuchi}
		H(\sigma_1)-H(\sigma)\leq 2U.
	\end{equation}
	Moreover, every triangular face in the internal boundary of the cluster shares one or two edges with the empty triangular faces outside the cluster, and we proceed as in the previous case. 
	\begin{itemize}
		\item If there exist two triangular faces $T, \, T'$ in the internal boundary of the cluster $C(\sigma)$ such that they both share two edges with some empty triangular faces, then we denote by $\eta$ the configuration obtained from $\sigma_1$ by exchanging the empty triangular face of the hole with $T$. We have
		\begin{align}
			& H(\eta)-H(\sigma_1)= -2U, \\
			& H(\eta)-H(\sigma)=[H(\eta)-H(\sigma_1)]+[H(\sigma_1)-H(\sigma)]\leq0.
		\end{align}
		Then, in a similar way, we move the empty triangular face $T_2$ until the internal boundary of the cluster in $T'$. During the first step, the energy possibly increases by $U$ as in the second case. If we denote by $\eta_1$ this configuration, then we have $H(\eta_1)-H(\eta)\leq U$. Moreover, $T'$ has two shared edges with some empty triangular face, then, denoting by $\xi$ the configuration obtained from $\eta_1$ by exchanging the empty triangular face of the hole with $T$, we have
		\begin{align}
			& H(\xi)-H(\eta_1)= -2U, \\
			& H(\xi)-H(\eta)=[H(\xi)-H(\eta_1)]+[H(\eta_1)-H(\eta)]\leq-U.
		\end{align}
		Thus we obtain
		\begin{equation}
			H(\xi)-H(\sigma)=[H(\xi)-H(\eta)]+[H(\eta)-H(\sigma)]=-U
		\end{equation}
		and therefore by \eqref{stab2U_nobuchi} it follows that $V_\sigma=2U$.
		\item If there exists only one triangular face $T$ such that it shares two edges with some triangular face outside of the cluster, then we define the configurations $\sigma_1$, $\eta$, $\eta_1$ and $\xi$ as before. Since now $H(\xi)-H(\eta_1)=-U$, we obtain $H(\xi)-H(\sigma)=0$. By construction $\xi$ has an internal angle of $\frac{5}{3}\pi$, see Figure \ref{anglesconf}(a). We define a configuration $\xi_1$ obtained from $\xi$ by getting in $\Lambda$ a new particle, and we define $\xi_2$ as the configuration obtained from $\xi_1$ by moving and attaching this particle to cove the internal angle of $\frac{5}{3}\pi$. We have
		\begin{equation}
			\begin{array}{ll}
				H(\xi_1)-H(\xi)&= \Delta, \\
				H(\xi_2)-H(\xi_1)& =-2U 
			\end{array}
		\end{equation}
		and therefore
		\begin{equation}
			H(\xi_2)-H(\sigma)=[H(\xi_2)-H(\xi_1)]+[H(\xi_1)-H(\xi)]+[H(\xi)-H(\sigma)]=-2U+\Delta.
		\end{equation}
		Thus by \eqref{stab2U_nobuchi} $V_\sigma=2U$.
		\item If each triangular face in the internal boundary has only one shared edge with an empty triangular face outside cluster, then we have 
		\begin{equation}
			H(\sigma_2)-H(\sigma_1)=-U,
		\end{equation} 
		where $\sigma_2$ is the configuration obtained from $\sigma_1$ by exchanging the empty triangular face of the hole with a triangular face $T_1$ in the internal boundary of the cluster. Thus, we obtain $H(\sigma_2)-H(\sigma)\leq U$, and by construction $\sigma_2$ has an internal angle of $\frac{5}{3}\pi$ (see Figure \ref{anglesconf}(a)). We define a configuration $\sigma_3$ obtained from $\sigma_2$ by getting in $\Lambda$ a new particle, and we define $\sigma_4$ as the configuration obtained from $\sigma_3$ by moving and attaching this particle to cover the internal angle of $\frac{5}{3}\pi$. We have
		\begin{equation}
			\begin{array}{ll}
				H(\sigma_3)-H(\sigma_2)&= \Delta, \\
				H(\sigma_4)-H(\sigma_3)&=-2U \\
			\end{array}
		\end{equation}
		and therefore
		\begin{equation}
        \begin{array}{ll}
			H(\sigma_4)-H(\sigma)&=[H(\sigma_4)-H(\sigma_3)]+[H(\sigma_3)-H(\sigma_2)]+[H(\sigma_2)-H(\sigma)] \\
            &=-U+\Delta>0.
        \end{array}
		\end{equation}
		We observe that in $\sigma_4$ there is an empty triangular face as in the second case. So, we iterate the same procedure starting from $\sigma_4$. Finally we obtain for the energy a total decreasing value $t\leq (-U+\Delta)+(-2U+\Delta)$, thus the stability level is $V_{\sigma}=\Delta+U$, which is obtained in $\sigma_3$.
	\end{itemize}
\end{proof*}

	\begin{figure}
	\centering
	\includegraphics[scale=0.5]{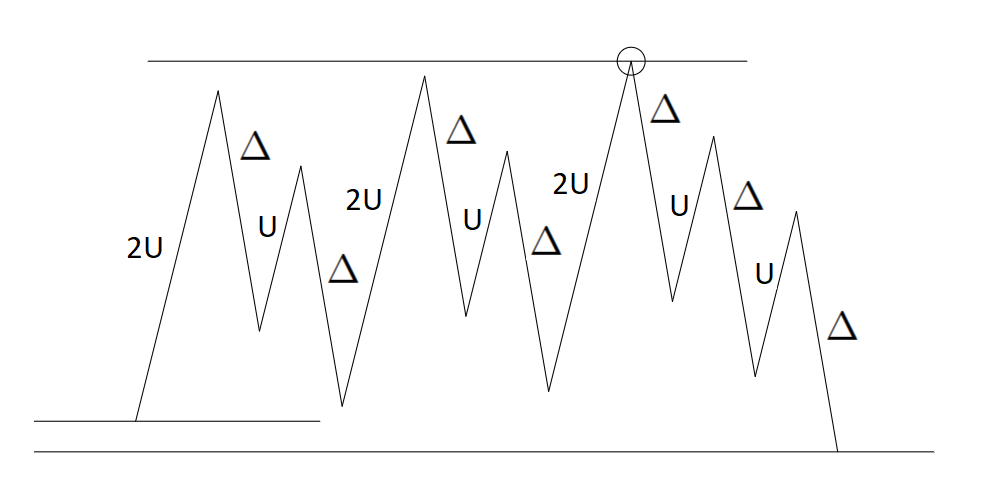}
	\caption{We depict an example of the energy landscape between $\cE_{B_m}(r)$ and $\cE_{B_{m-1}}(r)$ for the value $r=4$. We highlight with a circle the maximum of the energy, which is attained in $\o_9$.}
	\label{figure:Z_1}
\end{figure}

\begin{proof*}{\it  of Lemma \ref{zry}}
	We distinguish the three cases $\s\in Z$, $\s\in R$ and $\s\in Y$. Recall Definition \ref{corner} and extend it to clusters.  Moreover, recall the definition of $\Theta(\cdot,\cdot)$ given in Section \ref{defmodind} point 1.
	
	{\bf Stability level of $Z$.} 
	We begin by considering the set $Z_1$. For any configuration $\sigma\in Z_1$ we construct a path $\overline{\omega}\in \Theta(\sigma,I_{\sigma} \cap (Z_1 \cup \vuoto\}))$ that dismantles the bar on one of the shortest sides of the quasi-regular hexagon starting from one of its corners. Starting from $\sigma\equiv \omega_0\in Z_1$, we will define $\omega_1$ as follows. Consider a corner in one of the shortest sides of the cluster in $\mathcal{E}_{B_m}(r)$, with $m=0,...,5$ and let $j$ be a site belonging to this corner. Define $\omega_1$ as the configuration obtained starting from $\omega_0$ by moving the particle in $j$ to the empty nearest site. Since $\omega_1$ is obtained by breaking two bonds, we have $H(\omega_1)-H(\omega_0)=2U$. Then, consider $\omega_2$ as the configuration obtained from $\omega_1$ by moving the same particle outside $\Lambda$. We observe that $\omega_1$ and $\omega_2$ are not connected via one step of the dynamics, but there exist some configurations $\xi_1,...,\xi_n$ such that $(\omega_1,\xi_1,...,\xi_n,\omega_2)$ is a path with $H(\xi_i)=H(\xi_j)$ for all $i,j \in \{1,...,n\}$. We have $H(\omega_2)-H(\omega_1)=-\Delta$. Then, we analogously define $\omega_3$ and $\omega_4$ by considering the site $j_1$, where $j_1$ is the other site belonging to the same corner. In this case, when a particle is detached from the cluster defining $\omega_3$, only one bond is broken. Thus we have $H(\omega_3)-H(\omega_2)=U$ and $H(\omega_4)-H(\omega_3)=-\Delta$.
	By iterating this procedure along the considered side, a bar of the cluster is erased and we obtain the configuration $\eta \equiv \omega_k$ such that $\eta=\mathcal{E}_{B_{m-1}}(r)$ for $m \neq 0$, otherwise $\eta=\mathcal{E}_{B_5}(r-1)$ for $m=0$. Note that $k$ is twice the cardinality of the bar. See Figure \ref{figure:Z_1}.
	
	Note that if the initial configuration contains a regular hexagon with radius length one, then the final configuration contains a trapeze composed by three particles. 
	
	In order to determine where the maximum is attained, we observe that $H(\omega_{2j})<H(\omega_{2j+1})$ for every $j=0,...,\frac{k-2}{2}$ and $H(\omega_k)<H(\omega_{k-1})$. Thus, we will find the maximum over the configuration with odd index. By \eqref{def:hamiltonian}, we have
	\begin{align}
		& H(\omega_1)-H(\omega_0)=2U, \\
		& H(\omega_3)-H(\omega_0)=3U-\Delta
	\end{align}
	and for every $s=2,\cdots, \frac{k-4}{2}$, we have
	\begin{align}
		& H(\omega_{2s+1})-H(\omega_{2s-3})=3U-2\Delta, \\
		& H(\omega_{k-1})-H(\omega_{k-5})=2U-2\Delta.
	\end{align}
	It follows that for $n<k-1$ odd, $n=2\tilde s+1$, we obtain
	\begin{align}\label{differenza_hamiltoniana_Z1}
		H(\omega_n)-H(\omega_0) &=
		\begin{cases} 
			\displaystyle\sum_{\substack{s=2,..,\tilde s \\ s \, \text{even}}} [H(\omega_{2s+1})-H(\omega_{2s-3})] + [H(\omega_1)-H(\omega_0)] & \text{if } \tilde s \text{ even} \\ \\
			\displaystyle\sum_{\substack{s=3,..,\tilde s \\ s \, \text{odd}}} [H(\omega_{2s+1})-H(\omega_{2s-3})] + [H(\omega_3)-H(\omega_0)] & \text{if } \tilde s \text{ odd},
		\end{cases} \notag \\
		& = 
		\begin{cases} 
			\frac{\tilde s}{2}(3U-2\Delta) + 2U & \text{if } n=2\tilde s+1 \text{ with } \tilde s \text{ even},\\
			\frac{\tilde s -1}{2}(3U-2\Delta) + 3U-\Delta & \text{if } n=2\tilde s+1 \text{ with } \tilde s \text{ odd}.
		\end{cases} 
	\end{align}
	Thus, for $n=k-1$ we have
	\begin{align}\label{differenza_hamiltoniana_Z1_last}
		H(\omega_{k-1})-H(\omega_0) & =[H(\omega_{k-1})-H(\omega_{k-5})]+[H(\omega_{k-5})-H(\omega_0)] \notag \\
		& =2U-2\Delta+\frac{k-6}{4}(3U-2\Delta)+2U \notag \\
		& = \frac{k-2}{4}(3U-2\Delta)+U,
	\end{align}
	where we have used that $\tilde s$ is even, indeed $k-1=2(2r-j)-1=2 \tilde s-1$ with $j \in \{-1,+1,+3\}$.
	Since the result is an increasing function of $n=2\tilde s+1$, comparing the three maxima, we see that the absolute maximum is attained in $\omega_{k-5}$. Since $k$ is twice the cardinality of a bar, by Definitions \ref{def:bars} and \ref{qrhexagon}, we have
	\begin{itemize}
		\item $k=2(2r-1)$ if the initial configuration is $\mathcal{E}_{B_1}(r)$;
		\item $k=2(2r+1)$ if the initial configuration is $\mathcal{E}_{B_m}(r)$ for $m=2,3,4,5$;
		\item $k=2(2r+3)$ if the initial configuration is $\mathcal{E}(r+1)$.
	\end{itemize}
	So, by using \eqref{differenza_hamiltoniana_Z1} and replacing $k-5=2\tilde s +1$ with $\tilde s$ even, we have
	\begin{align}\label{phiomegabar}
		\Phi(\overline{\omega})-H(\omega_0) = H(\omega_{k-5})-H(\omega_0)=\frac{k-6}{4}(3U-2\Delta) + 2U.
	\end{align}
	Thus $\Phi(\overline{\omega})$ depends only on the value $k$, that is an increasing function of the radius $r$ of the quasi-regular hexagon. 
	The cardinality of the longest bar among those of the quasi-regular hexagon in a configuration in $Z_1$ is $2r^*+1$ (obtained by removing ${B_5}$ from $E_{B_5}(r^*)$), so we choose $k=2(2r^*+1)$. Note that the maximum is not obtained for $k=2(2r^*+3)$, since $\mathcal{E}(r^*+1) \not \in Z_1$. \\
	Let us check that $\omega_k \in I_{\sigma} \cap (Z_1 \cup \{\vuoto\})$. Since $k \leq 2(2r^*+1)$ with $r^*=\Big\lfloor \frac{U}{2(3U-2\Delta)} -1/2 \Big\rfloor$ and by \eqref{differenza_hamiltoniana_Z1_last}, we get
	\begin{align}\label{sottoI}
		H(\omega_0)-H(\omega_k) & =[H(\omega_0)-H(\omega_{k-1})]+[H(\omega_{k-1})-H(\omega_k)] \notag \\
		& =-\frac{k-2}{4}(3U-2\Delta)-U+\Delta \notag \\
		& \geq -\frac{2(2r^*+1)-2}{4}(3U-2\Delta)-U+\Delta =\d(3U-2\D)>0.
	\end{align}
	Finally, by equations \eqref{phiomegabar} and \eqref{sottoI}, we have
	\begin{equation}
		V_{\sigma} \leq \Phi(\overline{\omega}) -H(\sigma)=\frac{k-6}{4}(3U-2\Delta) + 2U.
	\end{equation}
	Thus, we find $V^*_{Z_1}=\max_{\sigma \in Z_1} V_{\sigma}$ by choosing $k=2(2r^*+1)$ and recalling $r^*=\Big\lfloor \frac{U}{2(3U-2\Delta)} -1/2 \Big\rfloor$, i.e.,
	\begin{equation}
		V^*_{Z_1} \leq 3\Delta-2U.
	\end{equation}
	
	Next, we analyze the set $Z_2$. For any configuration $\sigma\in Z_2$ we construct a path $\overline{\omega}\in \Theta(\sigma,I_{\sigma} \cap (Z_2 \cup \pieno\}))$. Starting from $\sigma\equiv \omega_0\in Z_2$, define $\omega_1$ by adding a free particle in $\Lambda$. Let us define $\omega_2$ in the following way. Consider a corner in one of the longest sides of the cluster in $\mathcal{E}_{B_m}(r)$ and let $j$ be a site belonging to this corner. Let $j_1$ be the site at distance one from $j$ such that $\sigma(j_1)=0$. We define $\omega_2$ by moving the free particle in $\omega_1$ until it reaches the site $j_1$. We observe that $\omega_1$ and $\omega_2$ are not connected via one step of the dynamics, but there exist some configurations $\xi_1,...,\xi_n$ such that $(\omega_1,\xi_1,...,\xi_n,\omega_2)$ is a path with $H(\xi_i)=H(\xi_j)$ for all $i,j \in \{1,...,n\}$. Moreover, we have
	\begin{align}
		& H(\omega_1)-H(\omega_0)= \Delta, \\
		& H(\omega_2)-H(\omega_1)=-U.
	\end{align}
	We consider $j_2$ the site at distance one from $j_1$ such that $\sigma(j_2)=0$ and $d(j_2,j')=2$ where $j' \neq j$ is another site of the initial cluster. The configuration $\omega_3$ is obtained from $\omega_2$ by adding a free particle, and $\omega_4$ is obtained from $\omega_3$ by moving the free particle until it reaches the site $j_2$. Again, we have 
	\begin{align}
		&H(\omega_3)-H(\omega_2)= \Delta, \\
		&H(\omega_4)-H(\omega_3)=-U.
	\end{align}
	Let us define $\omega_5$ and $\omega_6$. The configuration $\omega_5$ is obtained from $\omega_4$ by adding a free particle, and $\omega_6$ is obtained from $\omega_5$ by moving the free particle until it reaches the site $j_3$, where $j_3$ is the site at distance one from $j_2$ such that $\sigma(j_3)=0$ and $d(j_3,j')=1$ where $j' \neq j$ is another site of the initial cluster. We have
	\begin{align}
		& H(\omega_5)-H(\omega_4)= \Delta, \\
		& H(\omega_6)-H(\omega_5)=-2U.
	\end{align}
	We note that the energy has decreased by $2U$, since the particle has covered an internal angle of $\frac{5}{3}\pi$.
	By iterating this procedure along the considered side, a bar is added to the initial cluster. We obtain the configuration $\eta \equiv \omega_k$ such that $\eta= \mathcal{E}_{B_{m+1}}(r)$ for $m\neq 5$, otherwise $\eta=\mathcal{E}(r+1)$ for $m=5$. Note that the length of the path $k$ is equal to twice the cardinality of the bar. 
	
	In order to determine where the maximum is attained, we observe that $H(\omega_{2j})<H(\omega_{2j+1})$ for every $j=0,...,\frac{k-2}{2}$ and $H(\omega_k)<H(\omega_{k-1})$. Thus, we will find the maximum over the configuration with odd index. By \eqref{def:hamiltonian}, we have
	\begin{align}
		& H(\omega_1)-H(\omega_0)=\Delta, \\
		& H(\omega_3)-H(\omega_0)=2\Delta-U, \\
		& H(\omega_5)-H(\omega_0)=3\Delta-2U 
	\end{align}
	and for every $s=3,..., \frac{k-2}{2}$ we have
	\begin{equation}
		H(\omega_{2s+1})-H(\omega_{2s-3})=2\Delta-3U.
	\end{equation}
	It follows that for $n>5$ odd, $n=2\tilde s+1$, we obtain
	\begin{align}\label{differenza_hamiltoniana_Z2}
		H(\omega_n)-H(\omega_0) &=
		\begin{cases} 
			\displaystyle\sum_{\substack{s=3,..,\tilde s \\ s \, \text{odd}}} [H(\omega_{2s+1})-H(\omega_{2s-3})] + [H(\omega_5)-H(\omega_0)] & \text{if } \tilde s \text{ odd}, \notag \\ \\
			\displaystyle\sum_{\substack{s=4,..,\tilde s \\ s \, \text{even}}} [H(\omega_{2s+1})-H(\omega_{2s-3})] + [H(\omega_3)-H(\omega_0)] & \text{if } \tilde s \text{ even},
		\end{cases} \\
		& = 
		\begin{cases} 
			\frac{\tilde s-2}{2}(2\Delta-3U) + (3\Delta-2U) & \text{if } n=2\tilde s+1 \text{ with } \tilde s \text{ even},\\
			\frac{\tilde s -1}{2}(2\Delta-3U) + (2\Delta-U) & \text{if } n=2\tilde s+1 \text{ with } \tilde s \text{ odd}.
		\end{cases} 
	\end{align}
	Since $2\Delta-3U<0$ and therefore the result is a decreasing function of $n$, the absolute maximum is attained in $\omega_{5}$. So, we have
	\begin{align}\label{phiomegabarra}
		\Phi(\overline{\omega})-H(\omega_0)=H(\omega_5)-H(\omega_0)=3\Delta-2U.
	\end{align}
	Finally, let us check that $\omega_k \in I_{\sigma} \cap (Z_2 \cup \{\pieno\})$. If $\sigma \in Z_2 \setminus \mathcal{E}(r^*+1)$, then the cardinality of the smallest bar among those of the quasi-regular hexagon in a configuration in $Z_2$ is $k_{\text{min}}=2(r^*+1)+1$. 
	Since $r^*=\Big\lfloor \frac{U}{2(3U-2\Delta)} -\frac{1}{2} \Big\rfloor$ and by using \eqref{differenza_hamiltoniana_Z2}, we have 
	\begin{align}
		H(\omega_0)-H(\omega_k) & =[H(\omega_0)-H(\omega_{k-2})]+[H(\omega_{k-2})-H(\omega_k)]= \\
		&= \Big[\frac{k-1}{2}(2\Delta-3U) + (2\Delta-U)\Big]+U-\Delta \\
		&=(r^*+1)(2\Delta-3U) + (2\Delta-U)+U-\Delta= \\
		&=-2U+2\Delta>0,
	\end{align}
	since $\Delta>U$. Thus
	\begin{align}\label{Vsup}
		V_{\sigma} \leq \Phi(\overline{\omega})-H(\sigma)=3\Delta-2U.
	\end{align}
	Now we consider $\mathcal{E}(r^*+1)$ and we note that $H(\mathcal{E}(r^*+1))<H(\mathcal{E}_{B_1}(r^*+1))$. 
	Thus our path $\overline{\omega}$ is the composition of the path we have previously defined, which connects $\mathcal{E}(r^*+1)$ to $\mathcal{E}_{B_1}(r^*+1)$, and an additional part depending on the value of $\delta$ (recall that $\delta \in (0,1)$ is such that $r^*= \frac{U}{2(3U-2\Delta)} -\frac{1}{2}-\delta$). If $0<\delta<\frac{1}{2}$, then we add the bar $B_2$ as we have done above for $B_1$ obtaining that $\overline{\omega}$ connects $\mathcal{E}(r^*+1)$ to $\mathcal{E}_{B_2}(r^*+1)$ passing through $\mathcal{E}_{B_1}(r^*+1)$. If $\frac{1}{2}<\delta<1$, then in the same manner we add the bars $B_2,B_3,B_4,B_5,B_6$ obtaining that $\overline{\omega}$ conncets $\mathcal{E}(r^*+1)$ to $\mathcal{E}_{B_6}(r^*+1) \equiv \mathcal{E}(r^*+2)$ passing through $\mathcal{E}_{B_i}(r^*+1)$ for any $i=2,...,5$. In both cases the last configuration of the new paths belong to $\mathcal{I}_{\mathcal{E}(r^*+1)}$, indeed 
	\begin{align}
		H(\mathcal{E}(r^*+1))>H(\mathcal{E}_{B_2}(r^*+1)), & \qquad \text{if } \delta \in (0,\frac{1}{2}), \notag \\
		H(\mathcal{E}(r^*+1))>H(\mathcal{E}(r^*+2)), & \qquad \text{if } \delta \in (\frac{1}{2},1). \notag
	\end{align}
	Thus, using equations \eqref{differenza_hamiltoniana_Z2}, \eqref{phiomegabarra} and \eqref{Vsup}, we obtain
	
		\begin{align}
			V_\s&\leq
			\begin{cases}
				3\Delta-2U+ H(\mathcal{E}_{B_1}(r^*+1))-H(\mathcal{E}(r^*+1)), \; 
				\text{for } \delta \in \Big(0,\frac{1}{2}\Big), \\  
				3\Delta-2U+ H(\mathcal{E}_{B_5}(r^*+1))-H(\mathcal{E}(r^*+1)), \; 
				\text{for } \delta \in \Big(\frac{1}{2},1\Big),
			\end{cases} \notag \\
			&=
			\begin{cases}
				3\Delta-2U+\delta(3U-2\Delta)<2\Delta-\frac{U}{2}, \; 
				\text{for } \delta \in \Big(0,\frac{1}{2}\Big), \\  
				3\Delta-2U+(2\delta-1)(3U-2\Delta)<\Delta+U, \; 
				\text{for } \delta \in \Big(\frac{1}{2},1\Big). 
			\end{cases}
		\end{align}
	
	Thus we find
	\begin{equation}
		V^*_{Z_2}=\max_{\sigma \in Z_2} V_{\sigma} \leq \Delta+U.
	\end{equation}
	In conclusion, we have $V^*_Z=\max\{V^*_{Z_1},V^*_{Z_2}\} \leq \Delta+U$.
	
	{\bf Stability level of $R$.} Consider the set $R_1$. For any configuration $\sigma\in R_1$ we construct a path $\overline{\omega}\in \Theta(\sigma,I_{\sigma} \cap (R_1 \cup Z_1)))$. Starting from $\sigma\equiv \omega_0\in R_1$, let us define $\omega_1$ as follows. Consider the corner in one of the shortest sides of the cluster and let $j$ be a site belonging to it. Define the configuration $\omega_1$ starting from $\omega_0$ by moving the particle in $j$ to the nearest empty site. Since with this move two bonds are broken, we have $H(\omega_1)-H(\omega_0)=2U$. Then, consider $\omega_2$ as the configuration obtained from $\omega_1$ by moving the same particle outside $\Lambda$. We observe that $\omega_1$ and $\omega_2$ are not connected via one step of the dynamics, but there exist some configurations $\xi_1,...,\xi_n$ such that $(\omega_1,\xi_1,...,\xi_n,\omega_2)$ is a path with $H(\xi_i)=H(\xi_j)$ for all $i,j \in \{1,...,n\}$. We have $H(\omega_2)-H(\omega_1)=-\Delta$. 
	Then, analogously define $\omega_3$ and $\omega_4$ by considering the site $j'$, where $j'$ is the other site belonging to the same corner. In this case, when a particle is detached from the cluster defining $\omega_3$, only one bond is loss. Thus we have that $H(\omega_3)-H(\omega_2)=U$ and $H(\omega_4)-H(\omega_3)=-\Delta$.
	By iterating this procedure along the shortest side, a bar of the cluster is erased and we obtain the configuration $\eta \equiv \omega_l$, where $l$ is twice the cardinality of the considered bar. We observe that the greatest value of $l$ is always smaller than $k$, where $k$ is twice the cardinality of the greatest bar of the quasi-regular hexagon contained in the cluster.
	Analogously to the case $\o\in Z_1$, we have that $\omega_l \in I_{\sigma}$. Thus $V_{\sigma} < 3\Delta-2U$ and therefore
	\begin{equation}
		V^*_{R_1}= \max_{\sigma \in R_1} V_{\sigma} < 3\Delta-2U.
	\end{equation}
	Next, we consider the set $R_2$. For any configuration $\sigma\in R_2$ we construct a path $\overline{\omega}\in \Theta(\sigma,I_{\sigma} \cap (R_2 \cup Z_2 \cup \{\pieno \}))$. Starting from $\sigma\equiv \omega_0\in R_2$, let us define $\omega_1$ as follows. Consider a corner in one of the shortest sides of the cluster and let $j$ be a site belonging to it. We distinguish two cases depending on the length of the bar $l$ of the shortest side.
	
	\begin{itemize}
		\item If the cardinality of the bar $l$ is smaller than $2(r^*+1)-1$, we define $\omega_1$ by detaching the particle in $j$ from the cluster. Then, we define $\omega_2$ by moving the free particle outside $\Lambda$. Next, we consider the other site $j'$ belonging to the corner of the cluster and define $\omega_3$ and $\omega_4$ by detaching and moving the particle in $j'$ outside $\Lambda$. By iterating this procedure along the shortest side, a bar of the cluster is erased and we obtain the configuration $\eta \equiv \omega_{\tilde l}$, where $\tilde l=2l$ is twice the cardinality of the considered bar. Since $l< 2(r^*+1)-1$, we observe that the greatest value of $\tilde l$ is always smaller than $k$, where $k$ is twice the cardinality of the greatest bar of the quasi-regular hexagon contained in the cluster, that is $\tilde l < k$. Analogously to the case $\s\in Z_1$, we have that $\omega_{\tilde l} \in I_{\sigma}$. Thus $V_{\sigma} < 3\Delta-2U$.
		
		\item If the cardinality of the bar $l$ is at least $2(r^*+1)-1$, consider the site $j_1$ at distance one from $j$ and such that $\sigma(j_1)=0$. We define the configuration $\omega_1$ starting from $\omega_0$ and by adding a free particle. Then, we define $\omega_2$ by moving the free particle in $\omega_1$ until it reaches the site $j_1$. 
		Next, we consider $j_2$ the site at distance one from $j_1$ such that $\sigma(j_2)=0$ and $d(j_2,j')=2$, where $j' \neq j$ is another site of the initial cluster. We define $\omega_3$ and $\omega_4$ in the following way. The configuration $\omega_3$ is obtained from $\omega_2$ by adding a free particle in $\Lambda$, and $\omega_4$ is obtained from $\omega_3$ by moving the free particle until it reaches the site $j_2$. Let $j_3$ be the site at distance one from $j_2$ such that $\sigma(j_3)=0$ and $d(j_3,j')=1$, where $j' \neq j$ is another site of the initial cluster. We define $\omega_5$ and $\omega_6$ in the same way used before. The configuration $\omega_5$ is obtained from $\omega_4$ by adding a free particle in $\Lambda$, and $\omega_6$ is obtained from $\omega_5$ by moving the free particle until it reaches the site $j_3$.
		By iterating this procedure along the considered side, a bar is added to the initial cluster. Analogously to the case $\o\in Z_2$, we have that $\omega_{2l} \in I_{\sigma}$ since $l \geq 2(r^*+1)-1$. Thus $V_{\sigma} < \Delta+U$ and
		\begin{equation}
			V^*_{R_2} < \Delta+U.
		\end{equation}
	\end{itemize}
	In conclusion, we have that $V^*_R=\max\{V^*_{R_1},V^*_{R_2}\} < \Delta+U$. 
	
	{\bf Stability level of $Y$.}
	First, consider the set $Y_1$. For every configuration $\sigma$ in $Y_1$, all clusters are non-interacting and are of the same type of those in $Z_1$ or $R_1$. If $\sigma$ contains a cluster that is not a quasi-regular hexagon, then we take our path to be the path that cuts a bar, analogously to what has been done for $R_1$. We get a configuration in $\mathcal{I}_\sigma \cap Y_1$. Otherwise, if all clusters are quasi-regular hexagons, then we take our path to be the path that cuts a bar of the cluster, analogously to what has been done for $Z_1$. We get a configuration in $\mathcal{I}_\sigma \cap (Y_1 \cup Z_1)$. So, we have
	\begin{equation}
		V^*_{Y_1}=\max\{V^*_{R_1},V^*_{Z_1}\} <3\Delta-2U.
	\end{equation}
	Next, consider the set $Y_2$. For every configuration $\sigma$ in $Y_2$, there exists at least a cluster of the same type of those in $Z_2$ or $R_2$. If $\sigma$ contains a cluster of the type of those in $R_2$, i.e., $\sigma$ contains a cluster that is not a quasi-regular hexagon, we take our path to be the path that either cuts or adds a bar as it has been done for $R_2$. We get a configuration in $\mathcal{I}_\sigma \cap Y_2$. Otherwise, if the cluster is like those in $Z_2$, i.e., the cluster is a quasi-regular hexagon, then we take the path that adds a bar to the quasi-regular hexagon, alike the cases encountered when considering $Z_2$. We get a configuration in $\mathcal{I}_\sigma \cap (Y_2 \cup \{\pieno\})$. So, we have
	\begin{equation}
		V^*_{Y_2}=\max \{V^*_{R_2},V^*_{Z_2}\} < \Delta+U.
	\end{equation}
	We conclude that 
	\begin{equation*}
		V^*_Y=\max\{V^*_{Y_1},V^*_{Y_2}\}=V^*_Z.
	\end{equation*}
\end{proof*}

\subsection{Proof of Theorem \ref{thm:metastable_state}}
In this Section we identify stable and metastable states by proving Theorem \ref{thm:metastable_state}.

\begin{proof*}{\it  of Theorem~\ref{thm:metastable_state}}
	First, by direct computation we deduce that $H(\pieno)<H(\vuoto)$ if $L$ is sufficiently large, say $L>2r^*+3$. Moreover, we know that $\cX^s\subseteq\cX_V$ for any $V\geq0$. Thus, using Theorem \ref{prop:recurrence_property} and Proposition \ref{prop:stability_lower}, we conclude that $\cX^s=\{\pieno\}$. To show that $\cX^m=\{\vuoto\}$, we need to prove that $V_{\vuoto}=\Phi(\vuoto,\pieno)=\GK>V^*$, with $V^*=\D+U$. This part of the proof is analogue to the proof of \cite[equation (3.86)]{NOS}.
\end{proof*}

\subsection{Proof of Theorems \ref{thm:transition_time} and \ref{thm:gate}}\label{gateproof}
In this Section we give the proof of the main Theorems \ref{thm:transition_time} and \ref{thm:gate}.

\begin{proof*}{\it of Theorem~\ref{thm:transition_time}}
	Combining \cite[Theorem 4.1]{MNOS}, \cite[Theorem 4.9]{MNOS}, \cite[Theorem 4.15]{MNOS}, Theorem \ref{thm:metastable_state} and Corollary \ref{Gamma}, we get the claim.
\end{proof*}

\begin{proof*}{\it of Theorem~\ref{thm:gate}}
	If follows by Proposition \ref{prop:ingresso}.
\end{proof*}

\subsection{Proof of Theorem \ref{thm:subsup}}\label{subsup}
We refer to \cite[eq.\ (2.7)]{MNOS} for the definition of cycle. To prove Theorem \ref{thm:subsup} we need \cite[Theorem 3.2]{NZB}, which states that every state in a cycle is visited by the process before the exit with high probability. Using this result, to prove Theorem \ref{thm:subsup} we need to prove the following:
\bi
\item[1.] if $0<\d<\frac{1}{2}$, then 
\bi
\item[(i)] if $\h$ is a quasi-regular hexagon contained in $\cE_{B_4}(r^*)$, then there exists a cycle $\cC_{\vuoto}$ containing $\h$ and $\vuoto$ and not containing $\pieno$;
\item[(ii)] if $\h$ is a quasi-regular hexagon containing $\cE_{B_0}(r^*+1)$, then there exists a cycle $\cC_{\pieno}$ containing $\h$ and $\pieno$ and not containing $\vuoto$;
\ei
\item[2.]if $\frac{1}{2}<\d<1$, then 
\bi
\item[(i)] if $\h$ is a quasi-regular hexagon contained in $\cE_{B_0}(r^*+1)$, then there exists a cycle $\cC_{\vuoto}$ containing $\h$ and $\vuoto$ and not containing $\pieno$;
\item[(ii)] if $\h$ is a quasi-regular hexagon containing $\cE_{B_2}(r^*+1)$, then there exists a cycle $\cC_{\pieno}$ containing $\h$ and $\pieno$ and not containing $\vuoto$.
\ei
\ei

{\bf Case 1.} Let us start with (i). Let $\cC_{\vuoto}$ be the maximal connected set containing $\vuoto$ such that $\max_{\h'\in\cC_{\vuoto}}H(\h')<\GK$. Note that by definition $\cC_{\vuoto}$ is a cycle containing $\vuoto$ and not containing $\pieno$ since $\Phi(\vuoto,\pieno)=\GK$. It remains to prove that $\eta$ belongs to $\cC_{\vuoto}$. The proof goes as follows. We construct a path $\o^{\h,\vuoto}$ going from $\h$ to $\vuoto$ keeping the energy less than $\GK$. This path is obtained by erasing site by site each bar of $\h$, as explain in the first case of the proof of Lemma \ref{zry}. Let $\h\in\cE_{B_i}(r)$ with $r\leq r^*$ and $0\leq i\leq5$. If $\h\notin\cE_{B_0}(r)$, i.e., if $\h$ is not a regular hexagon, consider the sequence of configurations $\{\bar\o_i^{\h,\vuoto}\}_{i=-\bar{i},...,-1}$ connecting $\h$ to the regular hexagon $\cE_{B_0}(r)$ by erasing site by site each bar. If $\h\in\cE_{B_0}(r)$, we consider this path empty. From now on, let $\{\bar\o_i^{\h,\vuoto}\}_{i=0,...,r}$ be a sequence of configurations that contain regular hexagons, starting from $\cE_{B_0}(r)$ and ending in $\vuoto$, with radius $r-i$. To complete the construction we can use the same idea applied in the construction of the reference path. More precisely, between each pair $(\bar\o_i^{\h,\vuoto}, \bar\o_{i+1}^{\h,\vuoto})$ we can add a sequence of configurations $\tilde\o_i^{\h,\vuoto}=\{\tilde\o_{i,j}^{\h,\vuoto}\}_{j=0,...,12r-6}$ such that $\tilde\o_{i,0}^{\h,\vuoto}= \bar\o_i^{\h,\vuoto}$ and $\tilde\o_{i,j}^{\h,\vuoto}$ is obtained from $\bar\o_{i}^{\h,\vuoto}$ by erasing $j$ sites for $j>0$. Again, as in the reference path, the last interpolation consists in inserting between every pair of configurations in $\tilde\o_i^{\h,\vuoto}$ a sequence of configurations with a free particle in a suitable sequence of sites connecting the boundary of $\L$ to the site previously occupied by the erased particle. Either for any $r<r^*$ and $1\leq i\leq 6$ or $r=r^*$ and $1\leq i\leq 5$, we have that $H(\cE_{B_{i-1}}(r))<H(\cE_{B_5}(r^*))$. Thus by Proposition \ref{togliere} for the path $\o^{\h,\vuoto}$ we obtain 
\be{}
\max_{i}H(\o_i^{\h,\vuoto})=\max_{r\leq r^*}H(\cE_{B_{i-1}}(r))+3\D-2U<\GK.
\ee
The proof of (ii) is similar. Let $\cC_{\pieno}$ be the maximal connected set containing $\pieno$ such that $\max_{\h'\in\cC_\pieno}H(\h')<\GK$. Again $\cC_{\pieno}$ is a cycle containing $\pieno$ and not containing $\vuoto$ since $\Phi(\vuoto,\pieno)=\GK$. To prove that $\cC_\pieno$ contains $\h$ we define a path $\o^{\h,\pieno}$ going from $\h$ to $\pieno$ as follows. It is obtained first by reaching a regular hexagon shape and, from there, following the reference path $\o^*$ defined in Section \ref{refpath}. Suppose that $\h\in\cE_{B_i}(r)$, with $r\geq r^*$ and $0\leq i\leq5$. If $\h$ is a regular hexagon, then we define $\o^{\h,\pieno}$ as the part of the reference path going from $\h$ to $\pieno$. Otherwise, we add bars to $\h$ with a mechanism similar to the time reversal of the one used in the construction of $\o^{\h,\vuoto}$, until the path reaches a configuration in $\h\in\cE_{B_0}(r+1)$. The remaining part of the path follows the part of the reference path $\o^*$ from $\cE_{B_0}(r+1)$ to $\pieno$. Since for any $r\geq r^*$, $0\leq i\leq 5$ and $0<\d<\frac{1}{2}$, we have that $H(\cE_{B_i}(r))<H(\cE_{B_5}(r^*))$, for the path $\o^{\h,\pieno}$ we obtain
\be{}
\max_{i}H(\o_i^{\h,\pieno})=\max_{r\geq r^*}H(\cE_{B_i}(r))+3\D-2U<\GK.
\ee
{\bf Case 2.} The proof is analogue to the one done for the case 1 with the following changes. In the proof of (i), since $\frac{1}{2}<\d<1$, and either for any $r<r^*+1$ and $1\leq i\leq 6$ or for $r=r^*+1$ and $i=1$, we have that $H(\cE_{B_{i-1}}(r))<H(\cE_{B_1}(r^*+1))$. Thus by Proposition \ref{togliere} for the path $\o^{\h,\vuoto}$ we obtain
\be{}
\max_{i}H(\o_i^{\h,\vuoto})=\max_{r\leq r^*+1}H(\cE_{B_{i-1}}(r))+3\D-2U<\GK.
\ee
In the proof of (ii), since either for any $r>r^*+1$ and $0\leq i\leq 5$ or for $r=r^*+1$ and any $2\leq i \leq5$, we have that $H(\cE_{B_i}(r))<H(\cE_{B_1}(r^*+1))$, for the path $\o^{\h,\pieno}$ we obtain
\be{}
\max_{i}H(\o_i^{\h,\pieno})=\max_{r\geq r^*+1}H(\cE_{B_i}(r))+3\D-2U<\GK.
\ee
This concludes the proof.

\end{document}